\def\Resetstrings{
    \def\present{ }\let\bgroup={\let\egroup=}
    \def\Astr{}\def\astr{}\def\Atest{}\def\atest{}%
    \def\Bstr{}\def\bstr{}\def\Btest{}\def\btest{}%
    \def\Cstr{}\def\cstr{}\def\Ctest{}\def\ctest{}%
    \def\Dstr{}\def\dstr{}\def\Dtest{}\def\dtest{}%
    \def\Estr{}\def\estr{}\def\Etest{}\def\etest{}%
    \def\Fstr{}\def\fstr{}\def\Ftest{}\def\ftest{}%
    \def\Gstr{}\def\gstr{}\def\Gtest{}\def\gtest{}%
    \def\Hstr{}\def\hstr{}\def\Htest{}\def\htest{}%
    \def\Istr{}\def\istr{}\def\Itest{}\def\itest{}%
    \def\Jstr{}\def\jstr{}\def\Jtest{}\def\jtest{}%
    \def\Kstr{}\def\kstr{}\def\Ktest{}\def\ktest{}%
    \def\Lstr{}\def\lstr{}\def\Ltest{}\def\ltest{}%
    \def\Mstr{}\def\mstr{}\def\Mtest{}\def\mtest{}%
    \def\Nstr{}\def\nstr{}\def\Ntest{}\def\ntest{}%
    \def\Ostr{}\def\ostr{}\def\Otest{}\def\otest{}%
    \def\Pstr{}\def\pstr{}\def\Ptest{}\def\ptest{}%
    \def\Qstr{}\def\qstr{}\def\Qtest{}\def\qtest{}%
    \def\Rstr{}\def\rstr{}\def\Rtest{}\def\rtest{}%
    \def\Sstr{}\def\sstr{}\def\Stest{}\def\stest{}%
    \def\Tstr{}\def\tstr{}\def\Ttest{}\def\ttest{}%
    \def\Ustr{}\def\ustr{}\def\Utest{}\def\utest{}%
    \def\Vstr{}\def\vstr{}\def\Vtest{}\def\vtest{}%
    \def\Wstr{}\def\wstr{}\def\Wtest{}\def\wtest{}%
    \def\Xstr{}\def\xstr{}\def\Xtest{}\def\xtest{}%
    \def\Ystr{}\def\ystr{}\def\Ytest{}\def\ytest{}%
}
\Resetstrings

\def\Refformat{
         \if\Jtest\present
             {\if\Vtest\present\journalarticleformat
                  \else\conferencereportformat\fi}
            \else\if\Btest\present\bookarticleformat
               \else\if\Rtest\present\technicalreportformat
                  \else\if\Itest\present\bookformat
                     \else\otherformat\fi\fi\fi\fi}

\def\Rpunct{
   \def\Lspace{ }%
   \def\Lperiod{ }
   \def\Lcomma{ }
   \def\Lquest{ }
   \def\Lcolon{ }
   \def\Lscolon{ }
   \def\Lbang{ }
   \def\Lquote{ }
   \def\Lqquote{ }
   \def\Lrquote{ }
   \def\Rspace{}%
   \def\Rperiod{.}
   \def\Rcomma{,}
   \def\Rquest{?}
   \def\Rcolon{:}
   \def\Rscolon{;}
   \def\Rbang{!}
   \def\Rquote{'}
   \def\Rqquote{"}
   \def\Rrquote{`}
   }

\def\Lpunct{
   \def\Lspace{}%
   \def\Lperiod{\unskip.}
   \def\Lcomma{\unskip,}
   \def\Lquest{\unskip?}
   \def\Lcolon{\unskip:}
   \def\Lscolon{\unskip;}
   \def\Lbang{\unskip!}
   \def\Lquote{\unskip'}
   \def\Lqquote{\unskip"}
   \def\Lrquote{\unskip`}
   \def\Rspace{\spacefactor=1000}%
   \def\Rperiod{\spacefactor=3000}
   \def\Rcomma{\spacefactor=1250}
   \def\Rquest{\spacefactor=3000}
   \def\Rcolon{\spacefactor=2000}
   \def\Rscolon{\spacefactor=1250}
   \def\Rbang{\spacefactor=3000}
   \def\Rquote{\spacefactor=1000}
   \def\Rqquote{\spacefactor=1000}
   \def\Rrquote{\spacefactor=1000}
   }

\def\Refstd{
     \def\Acomma{\unskip, }
     \def\Aand{\unskip\ and }
     \def\Aandd{\unskip\ and }
     \def\Ecomma{\unskip, }
     \def\Eand{\unskip\ and }
     \def\Eandd{\unskip\ and }
     \def\acomma{\unskip, }
     \def\aand{\unskip\ and }
     \def\aandd{\unskip\ and }
     \def\ecomma{\unskip, }
     \def\eand{\unskip\ and }
     \def\eandd{\unskip\ and }
     \def\Namecomma{\unskip, }
     \def\Nameand{\unskip\ and }
     \def\Nameandd{\unskip\ and }
     \def\Revcomma{\unskip, }
     \def\Initper{.\ }
     \def\Initgap{\dimen0=\spaceskip\divide\dimen0 by 2\hskip-\dimen0}%
   }

\def\Smallcapsaand{
     \def\Aand{\unskip\bgroup{\Smallcapsfont\ AND }\egroup}%
     \def\Aandd{\unskip\bgroup{\Smallcapsfont\ AND }\egroup}%
     \def\eand{\unskip\bgroup\Smallcapsfont\ AND \egroup}%
     \def\eandd{\unskip\bgroup\Smallcapsfont\ AND \egroup}%
   }

\def\Smallcapseand{
     \def\Eand{\unskip\bgroup\Smallcapsfont\ AND \egroup}%
     \def\Eandd{\unskip\bgroup\Smallcapsfont\ AND \egroup}%
     \def\aand{\unskip\bgroup\Smallcapsfont\ AND \egroup}%
     \def\aandd{\unskip\bgroup\Smallcapsfont\ AND \egroup}%
   }

   \def\Citefont{}
   \def\ACitefont{}
   \def\Authfont{}
   \def\Titlefont{}
   \def\Tomefont{\sl}
   \def\Volfont{}
   \def\Flagfont{}
   \def\Reffont{\rm}
   \def\Smallcapsfont{\sevenrm}
   \def\Flagstyle#1{\hangindent\parindent\indent\hbox to0pt
       {\hss[{\Flagfont#1}]\kern.5em}\ignorespaces}

\def\Underlinemark{\vrule height .7pt depth 0pt width 3pc}

\def\Citebrackets{\Rpunct
   \def\Lcitemark{\def\Cfont{\Citefont}[\bgroup\Cfont}
   \def\Rcitemark{\egroup]}
   \def\LAcitemark{\def\Cfont{\ACitefont}\bgroup\ACitefont}%
   \def\RAcitemark{\egroup}
   \def\LIcitemark{\egroup}
   \def\RIcitemark{\bgroup\Cfont}
   \def\Citehyphen{\egroup--\bgroup\Cfont}
   \def\Citecomma{\egroup,\hskip0pt\bgroup\Cfont}%
   \def\Citebreak{}
   }

\def\Citeparen{\Rpunct
   \def\Lcitemark{\def\Cfont{\Citefont}(\bgroup\Cfont}
   \def\Rcitemark{\egroup)}
   \def\LAcitemark{\def\Cfont{\ACitefont}\bgroup\ACitefont}%
   \def\RAcitemark{\egroup}
   \def\LIcitemark{\egroup}
   \def\RIcitemark{\bgroup\Cfont}
   \def\Citehyphen{\egroup--\bgroup\Cfont}
   \def\Citecomma{\egroup,\hskip0pt\bgroup\Cfont}%
   \def\Citebreak{}
   }

\def\Citesuper{\Lpunct
   \def\Lcitemark{\def\Cfont{\Citefont}\raise1ex\hbox\bgroup\bgroup\Cfont}%
   \def\Rcitemark{\egroup\egroup}
   \def\LAcitemark{\def\Cfont{\ACitefont}\bgroup\ACitefont}%
   \def\RAcitemark{\egroup}
   \def\LIcitemark{\egroup\egroup}
   \def\RIcitemark{\raise1ex\hbox\bgroup\bgroup\Cfont}%
   \def\Citehyphen{\egroup--\bgroup\Cfont}
   \def\Citecomma{\egroup,\hskip0pt\bgroup%
      \Cfont}
   \def\Citebreak{}
   }

\def\Citenamedate{\Rpunct
   \def\Lcitemark{
      \def\Citebreak{\egroup\ [\bgroup\Citefont}
      \def\Citecomma{\egroup]; 
         \bgroup\let\uchyph=1\Citefont}(\bgroup\let\uchyph=1\Citefont}%
   \def\Rcitemark{\egroup])}
   \def\LAcitemark{
      \def\Citebreak{\egroup\ [\bgroup\Citefont}\def\Citecomma{\egroup], %
         \bgroup\ACitefont }\bgroup\let\uchyph=1\ACitefont}%
   \def\RAcitemark{\egroup]}
  \def\Citehyphen{\egroup--\bgroup\Citefont}
   \def\LIcitemark{\egroup}
   \def\RIcitemark{\bgroup\Citefont}
   }

\def\annotations{}
\def\Flagstyle#1{\par\noindent\hbox to 30pt{[#1]\hfil}%
\hangindent=30pt\hangafter=1}

\def\ifundefined#1{\expandafter\ifx\csname#1\endcsname\relax}
\ifundefined{thesis}\def\bye{\end{document}}\fi

\def\endauthor{:\kern6pt}

\Refstd\Citebrackets
\def\Citefont{\rm}\def\Titlefont{\Reffont}\def\Volfont{\bf}\def\Tomefont{\it}

\def\journalarticleformat{\Reffont\let\uchyph=1\parindent=1.25pc\def\Comma{}%

\sfcode`\.=1000\sfcode`\?=1000\sfcode`\!=1000\sfcode`\:=1000\sfcode`\;=1000\sfcode`\,=1000
                \par\vfil\penalty-200\vfilneg
      \if\Ftest\present\Flagstyle\Fstr\fi%

\if\Atest\present\bgroup\Authfont\Astr\egroup\def\Comma{\unskip\endauthor}\fi%
        \if\Ttest\present\Comma\bgroup\Titlefont\Tstr\egroup\def\Comma{, }\fi%
         \if\etest\present\hskip.2em(\bgroup\estr\egroup)\def\Comma{\unskip,
}\fi%
          \if\Jtest\present\Comma\bgroup\Tomefont\Jstr\/\egroup\def\Comma{,
}\fi%

\if\Vtest\present\if\Jtest\present\hskip.2em\else\Comma\fi\bgroup\Volfont\Vstr\egroup\def\Comma{, }\fi%
            \if\Dtest\present\hskip.2em(\bgroup\Dstr\egroup)\def\Comma{, }\fi%
             \if\Ptest\present\bgroup, \Pstr\egroup\def\Comma{, }\fi%
              \if\ttest\present\Comma\bgroup\Titlefont\tstr\egroup\def\Comma{,
}\fi%

\if\jtest\present\Comma\bgroup\Tomefont\jstr\/\egroup\def\Comma{, }\fi%

\if\vtest,\present\if\jtest\present\hskip.2em\else\Comma\fi\bgroup\Volfont\vstr\egroup\def\Comma{, }\fi%
                 \if\dtest\present\hskip.2em(\bgroup\dstr\egroup)\def\Comma{,
}\fi%
                  \if\ptest\present\bgroup, \pstr\egroup\def\Comma{, }\fi%
                   \if\Gtest\present{\Comma Gov't ordering no.
}\bgroup\Gstr\egroup\def\Comma{, }\fi%
                    \if\Mtest\present\Comma MR
\#\bgroup\Mstr\egroup\def\Comma{, }\fi%
                     \if\Otest\present{\Comma\bgroup\Ostr\egroup.}\else{.}\fi%
\annotations  
\vskip5ptplus1ptminus1pt}

\def\conferencereportformat{\Reffont\let\uchyph=1\parindent=1.25pc\def\Comma{}%

\sfcode`\.=1000\sfcode`\?=1000\sfcode`\!=1000\sfcode`\:=1000\sfcode`\;=1000\sfcode`\,=1000
                \par\vfil\penalty-200\vfilneg
      \if\Ftest\present\Flagstyle\Fstr\fi%

\if\Atest\present\bgroup\Authfont\Astr\egroup\def\Comma{\unskip\endauthor}\fi%
        \if\Ttest\present\Comma\bgroup\Titlefont\Tstr\egroup\def\Comma{, }\fi%
         \if\Jtest\present\Comma\bgroup\Tomefont\Jstr\/\egroup\def\Comma{,
}\fi%
          \if\Ctest\present\Comma\bgroup\Cstr\egroup\def\Comma{, }\fi%
           \if\Dtest\present\hskip.2em(\bgroup\Dstr\egroup)\def\Comma{, }\fi%
            \if\Mtest\present\Comma MR \#\bgroup\Mstr\egroup\def\Comma{, }\fi%
             \if\Otest\present{\Comma\bgroup\Ostr\egroup.}\else{.}\fi%
\annotations  
\vskip5ptplus1ptminus1pt}

\def\bookarticleformat{\Reffont\let\uchyph=1\parindent=1.25pc\def\Comma{}%

\sfcode`\.=1000\sfcode`\?=1000\sfcode`\!=1000\sfcode`\:=1000\sfcode`\;=1000\sfcode`\,=1000
                \par\vfil\penalty-200\vfilneg
      \if\Ftest\present\Flagstyle\Fstr\fi%

\if\Atest\present\bgroup\Authfont\Astr\egroup\def\Comma{\unskip\endauthor}\fi%
        \if\Ttest\present\Comma\bgroup\Titlefont\Tstr\egroup\def\Comma{, }\fi%
         \if\etest\present\hskip.2em(\bgroup\estr\egroup)\def\Comma{\unskip,
}\fi%
          \if\Btest\present\Comma in
\bgroup\Tomefont\Bstr\/\egroup\def\Comma{\unskip, }\fi%
           \if\otest\present\Comma\ \bgroup\ostr\egroup\def\Comma{, }\fi%
            \if\Etest\present\Comma\bgroup\Estr\egroup\unskip,
\ifnum\Ecnt>1eds.\else ed.\fi\def\Comma{, }\fi%
             \if\Stest\present\Comma\bgroup\Sstr\egroup\def\Comma{, }\fi%
              \if\Vtest\present\ \bgroup\bf\Vstr\egroup\def\Comma{, }\fi%
               \if\Ntest\present\Comma no. \bgroup\Nstr\egroup\def\Comma{,
}\fi%
                \if\Itest\present\Comma\bgroup\Istr\egroup\def\Comma{, }\fi%
                 \if\Ctest\present\ (\bgroup\Cstr\egroup)\def\Comma{, }\fi%
                  \if\Dtest\present\Comma\bgroup\Dstr\egroup\def\Comma{, }\fi%
                   \if\Ptest\present\Comma pp.\ \Pstr\def\Comma{, }\fi%

\if\ttest\present\Comma\bgroup\Titlefont\Tstr\egroup\def\Comma{, }\fi%
                     \if\btest\present\Comma in
\bgroup\Tomefont\bstr\egroup\def\Comma{, }\fi%
                       \if\atest\present\Comma\bgroup\astr\egroup\unskip,
\if\acnt\present eds.\else ed.\fi\def\Comma{, }\fi%
                        \if\stest\present\Comma\bgroup\sstr\egroup\def\Comma{,
}\fi%
                         \if\vtest\present\Comma vol.
\bgroup\vstr\egroup\def\Comma{, }\fi%
                          \if\ntest\present\Comma no.
\bgroup\nstr\egroup\def\Comma{, }\fi%

\if\itest\present\Comma\bgroup\istr\egroup\def\Comma{, }\fi%

\if\ctest\present\Comma\bgroup\cstr\egroup\def\Comma{, }\fi%

\if\dtest\present\Comma\bgroup\dstr\egroup\def\Comma{, }\fi%
                              \if\ptest\present\Comma\pstr\def\Comma{, }\fi%
                               \if\Gtest\present{\Comma Gov't ordering no.
}\bgroup\Gstr\egroup\def\Comma{, }\fi%
                                \if\Mtest\present\Comma MR
\#\bgroup\Mstr\egroup\def\Comma{, }\fi%

\if\Otest\present{\Comma\bgroup\Ostr\egroup.}\else{.}\fi%
\annotations  
\vskip5ptplus1ptminus1pt}

\def\bookformat{\Reffont\let\uchyph=1\parindent=1.25pc\def\Comma{}%

\sfcode`\.=1000\sfcode`\?=1000\sfcode`\!=1000\sfcode`\:=1000\sfcode`\;=1000\sfcode`\,=1000
                \par\vfil\penalty-200\vfilneg
      \if\Ftest\present\Flagstyle\Fstr\fi%

\if\Atest\present\bgroup\Authfont\Astr\egroup\def\Comma{\unskip\endauthor}%

\else\if\Etest\present\bgroup\def\Eand{\Aand}\def\Eandd{\Aandd}\Authfont\Estr\egroup\unskip, \ifnum\Ecnt>1eds.\else ed.\fi\def\Comma{, }%

\else\if\Itest\present\bgroup\Authfont\Istr\egroup\def\Comma{, }\fi\fi\fi%

\if\Ttest\present\Comma\bgroup\Tomefont\Tstr\/\egroup\def\Comma{\unskip, }%

\else\if\Btest\present\Comma\bgroup\Titlefont\Bstr\/\egroup\def\Comma{\unskip,
}\fi\fi%
            \if\otest\present\Comma\ \bgroup\ostr\egroup\def\Comma{, }\fi%

\if\etest\present\hskip.2em(\bgroup\estr\egroup)\def\Comma{\unskip, }\fi%
              \if\Stest\present\Comma\bgroup\Sstr\egroup\def\Comma{, }\fi%
               \if\Vtest\present\ \bgroup\bf\Vstr\egroup\def\Comma{, }\fi%
                \if\Ntest\present\Comma no. \bgroup\Nstr\egroup\def\Comma{,
}\fi%
                 \if\Atest\present\if\Itest\present
                         \Comma\bgroup\Istr\egroup\def\Comma{\unskip, }\fi%
                      \else\if\Etest\present\if\Itest\present
                              \Comma\bgroup\Istr\egroup\def\Comma{\unskip,
}\fi\fi\fi%
                     \if\Ctest\present\ (\bgroup\Cstr\egroup)\def\Comma{, }\fi%
                      \if\Dtest\present\Comma\bgroup\Dstr\egroup\def\Comma{,
}\fi%
                      \if\Ptest\present\bgroup, pp.\ \Pstr\egroup\def\Comma{,
}\fi

\if\ttest\present\Comma\bgroup\Titlefont\tstr\egroup\def\Comma{, }%

\else\if\btest\present\Comma\bgroup\Titlefont\bstr\egroup\def\Comma{, }\fi\fi%

\if\stest\present\Comma\bgroup\sstr\egroup\def\Comma{, }\fi%
                           \if\vtest\present\Comma vol.
\bgroup\vstr\egroup\def\Comma{, }\fi%
                            \if\ntest\present\Comma no.
\bgroup\nstr\egroup\def\Comma{, }\fi%

\if\itest\present\Comma\bgroup\istr\egroup\def\Comma{, }\fi%

\if\ctest\present\Comma\bgroup\cstr\egroup\def\Comma{, }\fi%

\if\dtest\present\Comma\bgroup\dstr\egroup\def\Comma{, }\fi%
                                \if\Gtest\present{\Comma Gov't ordering no.
}\bgroup\Gstr\egroup\def\Comma{, }\fi%
                                 \if\Mtest\present\Comma MR
\#\bgroup\Mstr\egroup\def\Comma{, }\fi%

\if\Otest\present{\Comma\bgroup\Ostr\egroup.}\else{.}\fi%
\annotations  
\vskip5ptplus1ptminus1pt}

\def\technicalreportformat{\Reffont\let\uchyph=1\parindent=1.25pc\def\Comma{}%

\sfcode`\.=1000\sfcode`\?=1000\sfcode`\!=1000\sfcode`\:=1000\sfcode`\;=1000\sfcode`\,=1000
                \par\vfil\penalty-200\vfilneg
      \if\Ftest\present\Flagstyle\Fstr\fi%

\if\Atest\present\bgroup\Authfont\Astr\egroup\def\Comma{\unskip\endauthor}%

\else\if\Etest\present\bgroup\def\Eand{\Aand}\def\Eandd{\Aandd}\Authfont\Estr\egroup\unskip, \ifnum\Ecnt>1eds.\else ed.\fi\def\Comma{, }%

\else\if\Itest\present\bgroup\Authfont\Istr\egroup\def\Comma{, }\fi\fi\fi%
          \if\Ttest\present\Comma\bgroup\Titlefont\Tstr\egroup\def\Comma{,
}\fi%
           \if\Atest\present\if\Itest\present
                   \Comma\bgroup\Istr\egroup\def\Comma{, }\fi%
                \else\if\Etest\present\if\Itest\present
                        \Comma\bgroup\Istr\egroup\def\Comma{, }\fi\fi\fi%
            \if\Rtest\present\Comma\bgroup\Rstr\egroup\def\Comma{, }\fi%
             \if\Ctest\present\ (\bgroup\Cstr\egroup)\def\Comma{, }\fi%
              \if\Dtest\present\Comma\bgroup\Dstr\egroup\def\Comma{, }\fi%
              \if\Ptest\present\bgroup, \Pstr\egroup\def\Comma{, }\fi%
               \if\ttest\present\Comma\bgroup\Titlefont\tstr\egroup\def\Comma{,
}\fi%
                \if\itest\present\Comma\bgroup\istr\egroup\def\Comma{, }\fi%
                 \if\rtest\present\Comma\bgroup\rstr\egroup\def\Comma{, }\fi%
                  \if\ctest\present\Comma\bgroup\cstr\egroup\def\Comma{, }\fi%
                   \if\dtest\present\Comma\bgroup\dstr\egroup\def\Comma{, }\fi%
                    \if\Gtest\present{\Comma Gov't ordering no.
}\bgroup\Gstr\egroup\def\Comma{, }\fi%
                     \if\Mtest\present\Comma MR
\#\bgroup\Mstr\egroup\def\Comma{, }\fi%
                      \if\Otest\present{\Comma\bgroup\Ostr\egroup.}\else{.}\fi%
\annotations  
\vskip5ptplus1ptminus1pt}

\def\otherformat{\Reffont\let\uchyph=1\parindent=1.25pc\def\Comma{}%

\sfcode`\.=1000\sfcode`\?=1000\sfcode`\!=1000\sfcode`\:=1000\sfcode`\;=1000\sfcode`\,=1000
                \par\vfil\penalty-200\vfilneg
      \if\Ftest\present\Flagstyle\Fstr\fi%

\if\Atest\present\bgroup\Authfont\Astr\egroup\def\Comma{\unskip\endauthor}%

\else\if\Etest\present\bgroup\def\Eand{\Aand}\def\Eandd{\Aandd}\Authfont\Estr\egroup\unskip, \ifnum\Ecnt>1eds.\else ed.\fi\def\Comma{, }%

\else\if\Itest\present\bgroup\Authfont\Istr\egroup\def\Comma{, }\fi\fi\fi%
          \if\Ttest\present\Comma\bgroup\Titlefont\Tstr\egroup\def\Comma{,
}\fi%
            \if\Atest\present\if\Itest\present
                    \Comma\bgroup\Istr\egroup\def\Comma{, }\fi%
                 \else\if\Etest\present\if\Itest\present
                         \Comma\bgroup\Istr\egroup\def\Comma{, }\fi\fi\fi%
                 \if\Ctest\present\ (\bgroup\Cstr\egroup)\def\Comma{, }\fi%
                  \if\Dtest\present\Comma\bgroup\Dstr\egroup\def\Comma{, }\fi%
                  \if\Ptest\present\bgroup, \Pstr\egroup\def\Comma{, }\fi%
                   \if\Gtest\present{\Comma Gov't ordering no.
}\bgroup\Gstr\egroup\def\Comma{, }\fi%
                    \if\Mtest\present\Comma MR
\#\bgroup\Mstr\egroup\def\Comma{, }\fi%
                     \if\Otest\present{\Comma\bgroup\Ostr\egroup.}\else{.}\fi%
\annotations  
\vskip5ptplus1ptminus1pt}

\catcode`\@=11
\def\matrixx#1{\null\,\vcenter{\normalbaselines\m@th
    \ialign{\hfil$\displaystyle ##$\hfil&&\quad\hfil$\displaystyle
##$\hfil\crcr
    \mathstrut\crcr\noalign{\kern-\baselineskip}
    #1\crcr\mathstrut\crcr\noalign{\kern-\baselineskip}}}\,}
\catcode`\@=12

\def\annotations{}
\def\Flagstyle#1{\par\noindent\hbox to 30pt{[#1]\hfil}%
\hangindent=30pt\hangafter=1}

\catcode`\@=11
\def\ATstartsection{\@startsection}
\def\resetATfont{\reset@font}
\def\zAT{\z@}
\def\RIfM@{\relax\protect\ifmmode}
\def\fraktur{\mathfrak}

\catcode`\@=12

\def\hexnumberZ#1{\ifnum#1<10 \number#1\else
 \ifnum#1=10 A\else\ifnum#1=11 B\else\ifnum#1=12 C\else
 \ifnum#1=13 D\else\ifnum#1=14 E\else\ifnum#1=15 F\fi\fi\fi\fi\fi\fi\fi}

\newcounter{alphactr}

\newcounter{romanctr}

\def\geom{g\'eom\'etrie}
\def\EXPOSE{expos\'e}
\def\SGA{S\'em\-in\-aire de G\'eom\'e\-trie Al\-g\'e\-bri\-que}
\def\PLUS{+}
\def\TIMES{*}

\def\begindiagram{\def\normalbaselines{\baselineskip20pt
      \lineskip3pt \lineskiplimit3pt}}

\catcode`\@=11
\def\operatoratfont{\operator@font}
\let\@afterindentfalse\@afterindenttrue
\@afterindenttrue
\catcode`\@=12

\def\op{\operatoratfont}

\def\manyxx#1{\if#1\PLUS{+}\else\if#1\TIMES{*}\fi\fi}

\catcode`\@=11
\def\backcases#1{\left.\vcenter{\normalbaselines\m@th
    \ialign{$##\hfil$\hfil\crcr#1\crcr}}\,\right\}}
\catcode`\@=12

\def\xhmode#1{\relax\ifmmode{#1}\else{\snap\hbox{$#1$}}\fi}
\def\snap{\hskip 0pt plus 4cm\penalty1000\hskip 0pt plus -4cm}
\def\PP{\xmode{\Bbb P\kern1pt}}
\def\email{jaffe{\kern0.5pt}@{\kern0.5pt}cpthree.unl.edu}
\def\N{\xmode{\Bbb N}}
\def\Z{\xmode{\Bbb Z}}
\def\O{{\cal O}}
\def\Gm{{\Bbb G}_m}
\def\Ga{{\Bbb G}_a}
\def\Coker{\mathop{\operatoratfont Coker}\nolimits}
\def\Ext{\mathop{\operatoratfont Ext}\nolimits}
\def\Hom{\mathop{\operatoratfont Hom}\nolimits}
\def\Ker{\mathop{\operatoratfont Ker}\nolimits}
\def\Mor{\mathop{\operatoratfont Mor}\nolimits}
\def\Nil{\mathop{\operatoratfont Nil}\nolimits}
\def\Pic{\mathop{\operatoratfont Pic}\nolimits}
\def\Spec{\mathop{\operatoratfont Spec}\nolimits}
\def\period{\mathop{\operatoratfont period}\nolimits}
\def\CHAR{\mathop{\operatoratfont char \kern1pt}\nolimits}
\def\SPEC{\mathop{\mathbf{Spec}}\nolimits}
\def\PIC{\mathop{\mathbf{Pic}}\nolimits}
\def\HOM{\mathop{\mathbf{Hom}}\nolimits}
\def\mp[[#1||#2||#3]]{\snap\hbox{$#1:\kern3pt #2\ \to\ #3$}}
\def\mapx[[#1||#2]]{\snap\hbox{$#1\ \to\ #2$}}
\def\dmap[[#1||#2||#3]]%
{     \setbox0=\hbox{$\displaystyle #1:\kern3pt #2\ \ \mapE\ \ \ #3$%
}     \ifdim\wd0<\hsize$$#1:\kern3pt #2\ \ \mapE\ \ \ #3$$\else%
     \begin{flushleft} $\displaystyle #1:\kern3pt #2$
     \end{flushleft} \begin{flushright} $\displaystyle\mapE{}\ \ \ #3$
     \end{flushright}\fi}
\def\dmapx[[#1||#2]]%
{     \setbox0=\hbox{$\displaystyle #1\ \ \mapE\ \ \ #2$%
}     \ifdim\wd0<\hsize$$#1\ \ \mapE\ \ \ #2$$\else%
     \begin{flushleft} $\displaystyle #1$ \end{flushleft}
     \begin{flushright} $\displaystyle\mapE{}\ \ \ #2$ \end{flushright}\fi}
\def\mapE#1{{\smash{\mathop{\longrightarrow}\limits^{#1}}}}
\def\smapE#1{{\smash{\mathop{\kern-10pt\rightarrow\kern-10pt}\limits^{#1}}}}
\def\mapS#1{\Big\downarrow\rlap{$\vcenter{\hbox{$\scriptstyle{#1}$}}$}}
\def\diagramx#1{{$$\begindiagram\matrixx{#1}$$}}
\def\diagramno#1#2{{$$\begindiagram\matrixx{#2}\eqno#1$$}}
\def\splitdiagram#1#2{
    \begin{flushleft}$\displaystyle\begindiagram\matrixx{#1}$\end{flushleft}
    \begin{flushright}$\displaystyle\begindiagram\matrixx{#2}$\end{flushright}}
\def\ses#1#2#3%
{  \setbox0=\hbox{$\matrixx{0&\mapE{}&#1&\mapE{}&#2&\mapE{}&#3&\mapE{}&0}$%
}  \ifdim\wd0<\hsize\diagramx{0&\mapE{}&#1&\mapE{}&#2&\mapE{}&#3&\mapE{}&
  0}\else\setbox0=\hbox{$\matrixx{0&\smapE{}&#1&\smapE{}&#2&\smapE{}&#3&
  \smapE{}&0}$}\ifdim\wd0<\hsize\diagramx{0&\smapE{}&#1&\smapE{}&#2&\smapE{}&
  #3&\smapE{}&0}\else\splitdiagram{0&\mapE{}&#1&\mapE{}&#2}{\mapE{}&#3
  &\mapE{}&0}\fi\fi}
\def\sesdot#1#2#3%
{  \setbox0=\hbox{$\matrixx{0&\mapE{}&#1&\mapE{}&#2&\mapE{}&#3&\mapE{}&0.}$%
}  \ifdim\wd0<\hsize\diagramx{0&\mapE{}&#1&\mapE{}&#2&\mapE{}&#3&\mapE{}&
  0.}\else\setbox0=\hbox{$\matrixx{0&\smapE{}&#1&\smapE{}&#2&\smapE{}&#3&
  \smapE{}&0.}$}\ifdim\wd0<\hsize\diagramx{0&\smapE{}&#1&\smapE{}&#2&\smapE{}&
  #3&\smapE{}&0.}\else\splitdiagram{0&\mapE{}&#1&\mapE{}&#2}{\mapE{}&#3
  &\mapE{}&0.}\fi\fi}
\def\sescomma#1#2#3%
{  \setbox0=\hbox{$\matrixx{0&\mapE{}&#1&\mapE{}&#2&\mapE{}&#3&\mapE{}&0,}$%
}  \ifdim\wd0<\hsize\diagramx{0&\mapE{}&#1&\mapE{}&#2&\mapE{}&#3&\mapE{}&
  0,}\else\setbox0=\hbox{$\matrixx{0&\smapE{}&#1&\smapE{}&#2&\smapE{}&#3&
  \smapE{}&0,}$}\ifdim\wd0<\hsize\diagramx{0&\smapE{}&#1&\smapE{}&#2&\smapE{}&
  #3&\smapE{}&0,}\else\splitdiagram{0&\mapE{}&#1&\mapE{}&#2}{\mapE{}&#3
  &\mapE{}&0,}\fi\fi}
\def\les#1#2#3%
{  \setbox0=\hbox{$\matrixx{0&\mapE{}&#1&\mapE{}&#2&\mapE{}&#3}$%
}  \ifdim\wd0<\hsize\diagramx{0&\mapE{}&#1&\mapE{}&#2&\mapE{}&
  #3}\else\setbox0=\hbox{$\matrixx{0&\smapE{}&#1&\smapE{}&#2&\smapE{}&
  #3}$}\ifdim\wd0<\hsize\diagramx{0&\smapE{}&#1&\smapE{}&#2&\smapE{}&
  #3}\else\splitdiagram{0&\mapE{}&#1&\mapE{}&#2}{\mapE{}&#3}\fi\fi}
\def\sesmaps#1#2#3#4#5{\setbox0=\hbox{$\matrixx{0&\mapE{}&#1&\mapE{#2}&#3&
     \mapE{#4}&#5&\mapE{}&0}$}\ifdim\wd0<\hsize\diagramx{0&\mapE{}&#1&
     \mapE{#2}&#3&\mapE{#4}&#5&\mapE{}&0}\else\setbox0=\hbox{$\matrixx{0&
     \smapE{}&#1&\smapE{#2}&#3&\smapE{#4}&#5&\smapE{}&
     0}$}\ifdim\wd0<\hsize\diagramx{0&\smapE{}&#1&\smapE{#2}&#3&\smapE{#4}&#5
     &\smapE{}&0}\else\splitdiagram{0&\mapE{}&#1&\mapE{#2}&#3}{\mapE{#4}&#5&
     \mapE{}&0}\fi\fi}
\def\sesmapsdot#1#2#3#4#5{\setbox0=\hbox{$\matrixx{0&\mapE{}&#1&\mapE{#2}&#3&
     \mapE{#4}&#5&\mapE{}&0.}$}\ifdim\wd0<\hsize\diagramx{0&\mapE{}&#1&
     \mapE{#2}&#3&\mapE{#4}&#5&\mapE{}&0.}\else\setbox0=\hbox{$\matrixx{0&
     \smapE{}&#1&\smapE{#2}&#3&\smapE{#4}&#5&\smapE{}&
     0.}$}\ifdim\wd0<\hsize\diagramx{0&\smapE{}&#1&\smapE{#2}&#3&\smapE{#4}&#5
     &\smapE{}&0.}\else\splitdiagram{0&\mapE{}&#1&\mapE{#2}&#3}{\mapE{#4}&#5&
     \mapE{}&0.}\fi\fi}
\def\seslab#1#2#3#4{$$0\ \ \mapE{}\ \ #1\ \ \mapE{}\ \ #2\ \ \mapE{}\ \ #3
                    \ \ \mapE{}\ \ 0\eqno(#4)$$}
\def\sesmapsone#1#2#3#4#5{$$1\ \ \mapE{}\ \ #1\ \ \mapE{#2}\ \ #3
                               \ \ \mapE{#4}\ \ #5
                               \ \ \mapE{}\ \ 1$$}
\def\seslabcomma#1#2#3#4{$$0\ \mapE{}\ #1\ \mapE{}\ #2\ \mapE{}\ #3
                    \ \mapE{}\ 0,\eqno(#4)$$}
\def\sesonerow#1#2#3{1&\mapE{}&#1&\mapE{}&#2&\mapE{}&#3&\mapE{}&1}
\def\sesonerowdot#1#2#3{1&\mapE{}&#1&\mapE{}&#2&\mapE{}&#3&\mapE{}&1\nulldot}
\def\makenull#1{{\setbox0=\hbox{#1}\wd0=0pt\box0}}
\def\nulldot{\makenull{.}}
\def\Rowfour#1#2#3#4{\diagramx{#1&\mapE{}&#2&\mapE{}&#3&\mapE{}&#4\cr}}
\def\rowfive#1#2#3#4#5{#1&\mapE{}&#2&\mapE{}&#3&\mapE{}&#4&\mapE{}&#5}
\def\Rowfive#1#2#3#4#5{\diagramx{#1&\mapE{}&#2&\mapE{}&#3&\mapE{}&#4&\mapE{}&#5\cr}}
\def\Rowsix#1#2#3#4#5#6{\diagramx{#1&\mapE{}&#2&\mapE{}&#3&\mapE{}&#4&\mapE{}&#5&\mapE{}&#6\cr}}
\def\Rowseven#1#2#3#4#5#6#7{\diagramx{#1&\mapE{}&#2&\mapE{}&#3&\mapE{}&#4&\mapE{}&#5&\mapE{}&#6&\mapE{}&#7\cr}}
\def\sRowsix#1#2#3#4#5#6{\diagramx{#1&\smapE{}&#2&\smapE{}&#3&\smapE{}&#4&\smapE{}&#5&\smapE{}&#6\cr}}
\def\o#1{\if#1\PLUS\oplus\else\if#1\TIMES\otimes\fi\fi}
\def\manyIN{\IN \cdots \IN}
\def\manytimes{\times \cdots \times}
\def\shC{{\cal{C}}}
\def\shM{{\cal{M}}}
\def\shN{{\cal{N}}}
\def\lfM{{\xmode{{\fraktur{\lowercase{M}}}}}}
\def\FG{finitely generated}
\def\WMAT{we may assume that}
\def\th#1{\xmode{{#1}^{\operatoratfont th}}}
\def\IN{\subset}
\def\NI{\supset}
\def\setof#1{\xmode{\{#1\}}}
\def\iso{\cong}
\def\qed{{\hfill$\square$}}
\def\RED#1{#1_{\operatoratfont red}}

\def\vec#1#2#3{#1_{#2}, \ldots, #1_{#3}}
\def\xmode#1{\relax\ifmmode{#1}\else{$#1$}\fi}
\def\cat#1{\xhmode{\ll\negthinspace\hbox{#1}\negthinspace\gg}}
\def\opcat#1{\xmode{\ll\negthinspace\hbox{#1}\negthinspace\gg^\circ}}
\def\pref#1{(\ref{#1})}
\def\block#1{\section{#1}}
\def\et#1{#1_{\hbox{\footnotesize\'et}}}
\def\Extfpqc{\Ext_{\op fpqc}}
\def\makeaddress{
     \vskip 0.15in
     \par\noindent {\footnotesize Department of Mathematics and Statistics,
                                  University of Nebraska}
     \par\noindent {\footnotesize Lincoln, NE 68588-0323, USA\ \ (\email)}}
\def\paidforby{ \def\thefootnote{\fnsymbol{footnote}}
     \par\noindent David B. Jaffe\protect\footnote{Partially supported by
                   the National Science Foundation.}
     \makeaddress\def\thefootnote{\arabic{footnote}}\setcounter{footnote}{0}}
\newenvironment{proof}{\trivlist \item[\hskip \labelsep{\sc
Proof.\kern1pt}]}{\endtrivlist}
 \newenvironment{proofnodot}{\trivlist \item[\hskip \labelsep{\sc
Proof}]}{\endtrivlist}
 \newenvironment{sketch}{\trivlist \item[\hskip \labelsep{\sc
Sketch.\kern1pt}]}{\endtrivlist}

\newenvironment{alphalist}{\begin{list}{(\alph{alphactr})}{\usecounter{alphactr}}}{\end{list}}

\newenvironment{romanlist}{\begin{list}{(\roman{romanctr})}{\usecounter{romanctr}}}{\end{list}}
 \newenvironment{definition}{\trivlist \item[\hskip
       \labelsep{\bf Definition.\kern1pt}]}{\endtrivlist}
\newenvironment{remarks}{\trivlist \item[\hskip
       \labelsep{\bf Remarks.\kern1pt}]}{\endtrivlist}
\newenvironment{problem}{\trivlist \item[\hskip
       \labelsep{\bf Problem.\kern1pt}]}{\endtrivlist}
\def\arrow(#1,#2){\ncline[nodesep=5pt]{->}{#1}{#2}}
\def\dottedarrow(#1,#2){\ncline[linestyle=dashed,nodesep=5pt]{->}{#1}{#2}}

\def\circno#1{{\bf [#1]}}
\def\smallcat#1{\cat{\kern3pt\fontsize{8}{10pt}\selectfont #1\kern3pt}}
\def\Morkschemes{\Mor_{\cat{\kern3pt\fontsize{8}{10pt}\selectfont
     $k$-schemes\kern3pt}}}
\documentclass[12pt]{article}\usepackage{amssymb}
\newtheorem{theorem}{Theorem}[section]
 \newtheorem{lemma}[theorem]{Lemma}
 \newtheorem{corollary}[theorem]{Corollary}
 \newtheorem{prop}[theorem]{Proposition}
 \newtheorem{exampleth}[theorem]{Example}
\newenvironment{example}{\begin{exampleth}\fontshape{n}\selectfont}{\end{exampleth}}

\begin{document}
\vskip 0.15in

\def\cabkf{\cat{(abelian group)-valued $k$-functors}}

\par\noindent{\Large\bf Functorial structure of units in a tensor product}
\vspace{0.15in}
\paidforby
\vspace{0.1in}

\block{Introduction}

Let $k$ be a field, and let $A$ be a \FG\ $k$-algebra.\footnote{All rings
in this paper are commutative.}  We explore the
structure of the functor from \cat{$k$-algebras} to \cat{abelian groups} given
by $B \mapsto (A \o*_k B)^*$.  More generally, if $S$ is a $k$-scheme of finite
type, not necessarily affine, we study the functor $\mu(S)$ given by
$B \mapsto (\Gamma(S,\O_S) \o*_k B)^*$.  This was done in
(\Lcitemark 8\Rcitemark \ 4.5) for the case where $k$ is algebraically
closed and $S$ is a variety.

We make the assumption that every irreducible component of $\RED{S}$ is
geometrically integral and has a rational point.  We summarize these
properties by saying that $S$ is {\it geometrically stable}.  If $S$ is
any $k$-scheme of finite type, we can always find a finite extension $k'$ of
$k$ such that $S \times_k \Spec(k')$ is geometrically stable as a $k'$-scheme.

With the assumption that $S$ is geometrically stable,
we find that $\mu(S)$ fits into an exact sequence
\ses{\Gm^r \times U \times \Z^n}{\mu(S)}{I%
}in which $I$ is a sheaf (for the fpqc topology), $I(B) = 0$ for every reduced
$k$-algebra $B$, and $U$ admits a finite filtration with successive quotients
isomorphic to $\Ga^{\kern1pt\beta}$, for various $\beta \in \N \cup
\setof{\infty}$.
We summarize these properties by saying that $I$ is {\it nilpotent\/} and
$U$ is {\it additive}.  In the sequence, $\Z^n$ denotes the constant sheaf
associated to the abelian group $\Z^n$, or equivalently, the functor which
represents the constant group scheme associated to the abelian group $\Z^n$.

Moreover, suppose we have a dominant morphism \hbox{\mp[[ f || S || T ]],} in
which both $S$ and $T$ are geometrically stable.  There is an
induced morphism of functors \mp[[ \mu(f) || \mu(T) || \mu(S) ]].  Let
$Q = \Coker[\mu(f)]$.  We find that $Q$ also fits into an exact sequence as
shown above, except that $\Gm^r \times U \times \Z^n$ is replaced by an
extension of a \FG\ abelian group (i.e.\ the associated constant sheaf) by
$\Gm^r \times U$, $U$ is pseudoadditive (see p.\ \pageref{pseudoadditive-def}),
and we do not know if $I$ is a sheaf.  Correspondingly,
we do not know if $Q$ is a sheaf, but we do know at least that
$Q|_{\smallcat{reduced $k$-algebras}}$ is a sheaf and moreover that the
canonical map \mapx[[ Q || Q^+ ]] is a monomorphism.

Specializing to the affine case, we see for example that if $A$ is a subalgebra
of a $k$-algebra $C$ (and $\Spec(A)$, $\Spec(C)$ are geometrically stable),
then the functor given by $B \mapsto (C \o* B)^*/(A \o* B)^*$ fits into such an
exact sequence.

We have thus far described the content of the first theorem
\pref{tori-result-generalized} of this paper.  Now we describe the second
theorem \pref{kernel-pic-nilimmersion}, which is an application of the first.

Let $X$ be a geometrically stable $k$-scheme.  Let
\hbox{\mp[[ i || X_0 || X ]]} be a nilimmersion, such that the ideal sheaf
$\shN$ of $X_0$ in $X$ has square zero.  Let $P$ be the functor from
\cat{$k$-algebras} to \cat{abelian groups} given by
$$P(B) = \Ker[ \Pic(X \times_k \Spec(B))\ \rightarrow
         \ \Pic(X_0 \times_k \Spec(B))].$%
$Of course, if $X$ is affine, $P = 0$, but in general $P$ is not zero.  We
find that $P$ fits into an exact sequence
\sescomma{D \o+ I}{U}{P%
}in which $I$ is nilpotent (except possibly not a sheaf), $U$ is
pseudoadditive,
and $D$ is the constant sheaf associated to a \FG\ abelian group.

Although this theorem does not imply that $P$ is a sheaf, it does imply that
if \mp[[ f || B || C ]] is a faithfully flat homomorphism of reduced
$k$-algebras, then $P(f)$ is injective \pref{sheaf-kernel-pic-nilimmersion}.
In fact, this holds even if $\shN^2 \not= 0$.

We indicate the idea of the proof of the second theorem.  We have an exact
sequence
\splitdiagram{H^0(X,\O_X^*)&\mapE{}&H^0(X_0,\O_{X_0}^*)%
}{\mapE{}&H^1(X,\shN)&\mapE{}&
\Ker[\Pic(X)\ \rightarrow\ \Pic(X_0)]&\mapE{}&0.%
}Functorializing this yields an exact sequence:
\sescomma{\Coker[\mu(i)]}{\Ga^{\kern1pt\beta}}{P%
}in which $\beta = h^1(X,\shN)$.  The first theorem tells us what
$\Coker[\mu(i)]$ is like.  The second theorem is deduced from this.

Finally, we describe a theorem about the Picard group, whose proof in
\Lcitemark 7\Rcitemark \Rspace{} uses both theorems of this paper.  Let $k$
be a field, and let
$X$ be a separated $k$-scheme of finite type.  Then there exists a finite
field extension $k^+$ of $k$ such that for every algebraic extension $L$ of
$k^+$, the canonical map \mapx[[ \Pic(X_L) || \Pic(X_{L^a}) ]] is injective.

\vspace{0.1in}

\par\noindent{{\bf Acknowledgements.}\ Bob Guralnick supplied the
neat proof of \pref{unit-lemma-generalized}.  Faltings kindly provided example
\pref{Faltings}, thereby correcting an error.

\vspace{0.1in}

\par\noindent{\bf Conventions.}
\begin{alphalist}
\item A {\it $k$-functor\/} is a functor from \cat{$k$-algebras} to \cat{sets}.
(The usage of the term {\it $k$-functor\/} here is slightly different from the
usage in\Lspace \Lcitemark 8\Rcitemark \Rspace{}.) If $V$ is a $k$-scheme, then
we also
let $V$ denote the representable $k$-functor given by
$V(B) = \Morkschemes(\Spec(B),V)$.
\item A $k$-functor $F$ is a {\it sheaf\/} (by which we mean
{\it sheaf for the fpqc topology}) if for every faithfully flat homomorphism
\mp[[ p || B || C ]], the canonical map\label{Psi-place}
\dmap[[ \Psi_{F,p} || F(B) || \setof{x \in F(C): F(i_1)(x) = F(i_2)(x)} ]]%
is bijective, where \mp[[ i_2, i_2 || C || C \o*_B C ]] are given by
$c \mapsto c \o* 1$ and $c \mapsto 1 \o* c$, respectively.
\item The superscript $+$ is used to denote {\it associated sheaf}.
\item If $k$ is a field, $X$ is a $k$-scheme, and $L$ is a field extension of
      $k$, we let $X_L$ denote $X \times_k \Spec(L)$.  We let $k^a$ denote an
      algebraic closure of $k$.
\item If $X$ is a scheme, we let $\Gamma(X)$ denote $\Gamma(X,\O_X)$, and we
      let $\Gamma^*(X)$ denote $\Gamma(X,\O_X)^*$.
\item If $B$ is a ring, $\Nil(B)$ denotes its nilradical.
\item $k$-functors are said to be
      {\it (abelian group)-valued\/} if they take values in
      \cat{abelian groups} rather than \cat{sets}.
\end{alphalist}

We give some definitions which are adapted from
\Lcitemark 8\Rcitemark \Rspace{}\ pp.\ 173, 180.  If $k$ is a field and $X$,
$Y$
are $k$-schemes, then $\HOM(X,Y)$ denotes the $k$-functor given by
$$B \mapsto \Morkschemes(X \times_k \Spec(B), Y).$%
$
An (abelian group)-valued $k$-functor $F$ is {\it nilpotent\/} if it is a sheaf
and $F(B) = 0$ for every reduced $k$-algebra $B$.  We say that $F$ is
{\it subnilpotent\/} if it can be embedded as a subsheaf of a nilpotent
$k$-functor.

An (abelian group)-valued $k$-functor is {\it discrete\/} if it is a constant
sheaf.  We also say that such a functor is (for example)
{\it discrete and \FG}, if it is the constant sheaf associated to a
\FG\ abelian group.

\block{Additive $k$-functors}

Let $k$ be a field.  We need to consider a countably-infinite-dimensional
analog of unipotent group schemes over $k$.  Actually, what we will be
considering is more restrictive, as we shall only be considering the analog of
unipotent group schemes over $k$ which are smooth, connected, and moreover
which are $k$-solvable.  (See\Lspace \Lcitemark 10\Rcitemark \Rspace{}\ \S5.1.)

The simplest infinite dimensional example is the (abelian group)-valued
$k$-functor $\Ga^\infty$ given by $B \mapsto \o+_{i=1}^\infty B$.  However,
this is not good enough for our purposes, since in positive characteristic
one can have nontrivial extensions of $\Ga$ by $\Ga$.  (See
e.g.{\ }\Lcitemark 13\Rcitemark \ p.\ 67, exercise 8 or
\Lcitemark 12\Rcitemark \Rspace{}\ VII\ \S2.)
We want to define a class of objects which is closed under extension.

\begin{definition}
Let $k$ be a field, and let $F$ be an (abelian group)-valued $k$-functor.  Then
$F$ is {\it strictly additive\/} if it is isomorphic to $\Ga^\alpha$ for some
$\alpha \in \setof{0,1,\ldots,\infty}$, and $F$ is {\it additive\/} if it
admits a filtration:
$$0\ =\ F_0\ \IN F_1\ \manyIN F_n\ =\ F,$%
$whose successive quotients are strictly additive.
\end{definition}

This terminology is not perfect, but it is at least consistent with the
usage of the word {\it additive\/} in\Lspace \Lcitemark 8\Rcitemark \Rspace{}:
by (\ref{ext-results}\ref{additive-additive-char0}), it will follow that if
$k$ has characteristic zero, then additive $\Longrightarrow$ strictly additive.

We define the {\it dimension\/} of an additive $k$-functor $F$ to be the sum
of the dimensions of the successive quotients in a filtration of $F$, as
in the definition of additive.  Thus $\dim(F) \in \setof{0,1,\ldots,\infty}$.

If $F$ is additive, we define its {\it period\/} to be the smallest $n$ for
which there exists a filtration as in the definition of additive.

A direct sum of countably many additive $k$-functors need not be additive,
even if the summands are finite-dimensional.  Also, we shall not concern
ourselves with uncountable direct sums (e.g.\ of $\Ga$), as they seem not
to arise in practice.  The following two statements are easily checked:

\begin{prop}\label{extension-of-additive}
Let $k$ be a field.  Let
\ses{F'}{F}{F''%
}be an exact sequence of (abelian group)-valued $k$-functors, in which
$F'$ and $F''$ are additive.  Then $F$ is additive.
\end{prop}

\begin{prop}\label{quotient-of-additive}
Let $k$ be a field of characteristic zero.  Let
\ses{F'}{F}{F''%
}be an exact sequence of (abelian group)-valued $k$-functors, in which
$F'$, $F$ are additive.  Then $F''$ is additive.
\end{prop}

\begin{prop}\label{lemma-two}
Let $k$ be a field.  Let
\ses{F'}{F}{F''%
}be an exact sequence of (abelian group)-valued $k$-functors, in which
$F'$, $F$ are additive and finite-dimensional.  Then $F''$ is additive.
\end{prop}

\begin{proof}
By (\Lcitemark 3\Rcitemark \ 11.17), $(F'')^+$ is representable.  Let
\mp[[ p || F || (F'')^+ ]] be the canonical map, which is fpqc-surjective.
By (\Lcitemark 11\Rcitemark \ Theorem 10), there exists a morphism
\mp[[ \sigma || (F'')^+ || F ]] of $k$-functors such that
$p \circ \sigma = 1_{(F'')^+}$.  Hence $F'' = (F'')^+$.  Let $X'$, $X$, and
$X''$ be the group schemes which represent $F'$, $F$, and $F''$, respectively.
Then we have an exact sequence
\ses{X'}{X}{X''%
}in \cat{commutative $k$-group schemes}.  Since $F$ is additive and
finite-dimensional, $X$ admits a series whose factors are copies of the
group scheme $\Ga$.  Hence $X''$ admits a series whose factors are group scheme
quotients of $\Ga$.  But any quotient of $\Ga$ is $0$ or $\Ga$
(\Lcitemark 10\Rcitemark \ 2.3), so $X''$ admits a series whose
factors are the group scheme $\Ga$.  From the argument at the beginning of
the proof (showing that under certain circumstances
fpqc-surjective $\Longrightarrow$ surjective), we see that $F''$ admits a
series whose factors are the (abelian group)-valued $k$-functor $\Ga$.
Hence $F''$ is additive.  \qed
\end{proof}

Unfortunately, \pref{quotient-of-additive} fails in positive characteristic.
We will give an example of this, but there are a couple of preliminaries:

\begin{lemma}\label{strictly-additive-splitting}
If $F$ and $G$ are strictly additive and \mp[[ \pi || F || G ]] is an
epimorphism in \cabkf, then $\pi$ splits.
\end{lemma}

\begin{sketch}
The lemma is clear if $\CHAR(k) = 0$, so \WMAT\ $\CHAR(k) = p > 0$.
We do the case where $F = G = \Ga^\infty$; the proof in the other cases
is the same.  Let $e_i$ denote the element of $\Ga^\infty(k)$
which has $1$ in the \th{i} spot and $0$'s elsewhere.  Let $A = k[t]$.
For each $i$, choose $f_i \in \Ga^\infty(A)$ such that $\pi(f_i) = te_i$.
Let $g_i$ be the part of $f_i$ involving only the monomials
$t,t^p,t^{p^2},\ldots$.  Since the monomials which appear in an expression
for $\pi$ also have this form (with various $t$), it follows that
$\pi(g_i) = te_i$.

Regard $g_i$ as a function of $t$.  Define \mp[[ \sigma || G || F ]] by
$\sigma(be_i) = g_i(b)$, where $B$ is a $k$-algebra and $b \in B$.  Then
$\sigma$ splits $\pi$.  \qed
\end{sketch}

\begin{corollary}\label{strictly-additive-quotient}
If $F$ is strictly additive and $G$ is additive, and
\hbox{\mp[[ \pi || F || G ]]}
is an epimorphism in \cabkf, then $G$ is strictly additive.
\end{corollary}

\begin{proof}
Induct on the period of $G$.  If $\period(G) \leq 1$ we are done.  Otherwise,
we can find $G' \IN G$ and an exact sequence
\sesmaps{G'}{}{G}{p}{H%
}in which $H$ is strictly additive, $G'$ is additive, and
$\period(G') < \period(G)$.  By \pref{strictly-additive-splitting},
$p \circ \pi$ splits.  Hence $p$ splits.  Hence $\period(G) = \period(G')$:
contradiction.  \qed
\end{proof}

\begin{example}\label{Faltings}
(provided by G.\ Faltings)
\par\noindent Let \mp[[ f || \Ga^\infty || \Ga^\infty ]] be given by
$(x_1,x_2,x_3,\ldots) \mapsto (x_1, x_2 - x_1^p, x_3 - x_2^p, \ldots)$.
Then $f$ is a monomorphism.  Let $F'' = \Coker(f)$.   If $F''$ were additive,
then by \pref{strictly-additive-quotient} $F''$ would be strictly additive,
and so by \pref{strictly-additive-splitting} $f$ would split.  However, this
is clearly not the case.  Hence $F''$ is not additive.
\end{example}

The following generalization of {\it additive\/} allows us to work around the
behavior illustrated by the example:

\begin{definition}\label{pseudoadditive-def}
An (abelian group)-valued $k$-functor $P$ is {\it pseudoadditive\/} if
for some $n \in \N$ there exists an exact sequence
\Rowseven{0}{U_1}{U_2}{\cdots}{U_n}{P}{0%
}in \cabkf\ in which $\vec U1n$ are additive.
\end{definition}

By example \pref{Faltings}, one cannot always take $n=1$, i.e.\ additive
$\not=$ pseudoadditive.  We do not know if one can always take $n=2$.

\begin{prop}\label{pseudoadditive-is-sheaf}
If $P$ is pseudoadditive, then $P$ is a sheaf.
\end{prop}

\begin{proof}
Let $\shC$ be the class of (abelian group)-valued $k$-functors $F$ with the
property that for any faithfully flat homomorphism \mapx[[ B || C ]] of
$k$-algebras, the usual \v Cech complex
\sRowsix{0}{F(B)}{F(C)}{F(C \o*_B C)}{F(C \o*_B C \o*_B C)}{\cdots%
}is exact.  Then $\Ga^\alpha \in \shC$ for all $\alpha$.  If
\ses{F'}{F}{F''%
}is an exact sequence and any two of $F', F, F''$ are in $\shC$, then so is
the third.  Hence $P$ is in $\shC$, so $P$ is a sheaf.  \qed
\end{proof}

\block{Extensions in \cabkf}

For any objects $F_1, F_2$ in an abelian category, one can define an abelian
group $\Ext^1(F_1, F_2)$, whose elements are isomorphism classes of extensions
\sesdot{F_2}{F}{F_1%
}(The general theory is described in
\Lcitemark 9\Rcitemark \Rspace{}\ Ch.\ VII, among other places.)  Also, we will
refer
to such an exact sequence as defining an {\it extension of $F_1$ by $F_2$}.

In particular, the theory applies to \cabkf.
We shall say that an exact sequence
\sesmaps{F_2}{}{F}{\pi}{F_1%
}in this category is {\it set-theoretically split\/} if there exists a morphism
of $k$-functors \mp[[ \sigma || F_1 || F ]] such that
$\pi \circ \sigma = 1_{F_1}$.  We also refer to
{\it set-theoretically split extensions}.  Let $\Ext^1_s(F_1,F_2)$ denote
the subgroup of elements of $\Ext^1(F_1,F_2)$ which correspond to
set-theoretically split extensions.

To compute $\Ext^1_s(F_1,F_2)$, we copy (with appropriate but minor changes)
some definitions
which may be found in (\Lcitemark 12\Rcitemark \ VII\ \S4).
For this discussion, fix $F_1$ and $F_2$.  A {\it symmetric factor system\/} is
a morphism \mp[[ f || F_1 \times F_1 || F_2 ]] of $k$-functors such that
\begin{eqnarray*}
0 & = & f(y,z) - f(x+y,z) + f(x,y+z) - f(x,y)\\
f(x,y) & = & f(y,x)
\end{eqnarray*}
for all $k$-algebras $B$ and all $x,y,z \in F_1(B)$.  If
\mp[[ g || F_1 || F_2 ]] is a morphism of $k$-functors, then there is a
symmetric factor system $\delta g$ defined by
$$\delta g(x,y) = g(x+y) - g(x) - g(y);$%
$such a system is called {\it trivial}.  The group structure on $F_2$ makes
the set of symmetric factor systems into a group.  Then by standard arguments,
$\Ext^1_s(F_1,F_2)$ is isomorphic to the group of symmetric factor systems,
modulo the subgroup of trivial factor systems.

If $F_1$ and $F_2$ are sheaves, then one can also compute the group
$\Ext^1_{\op fpqc}(F_1,F_2)$, i.e.\ the group of isomorphism classes of
extensions of $F_1$ by $F_2$ in
\cat{(abelian group)-valued $k$-functors which are sheaves}.  If moreover
$F_1$ and $F_2$ are represented by commutative group schemes $X_1$ and $X_2$
of finite type over $k$, then [see\Lspace \Lcitemark 4\Rcitemark \Rspace{}\ 5.4
and
\Lcitemark 2\Rcitemark \Rspace{}\ 3.5, 7.3(ii)]
$\Ext^1_{\op fpqc}(F_1,F_2) = \Ext^1(X_1,X_2)$, where
the latter $\Ext$ group is computed relative to the abelian category
\cat{commutative group schemes of finite type over $k$}.  For arbitrary
sheaves $F_1$, $F_2$, we have $\Ext^1(F_1,F_2) \IN \Ext^1_{\op fpqc}(F_1,F_2)$,
but not equality in general, as may be seen e.g.\ from the exact sequence
\sesmaps{\Z/p\Z}{}{\Ga}{t\ \mapsto\ t^p - t}{\Ga%
}in the group scheme category, where $\CHAR(k) = p > 0$.

\begin{prop}\label{ext-results}
Let $k$ be a field.  We consider objects and morphisms in
\cabkf.  Then:
\begin{alphalist}
\item If \mapx[[ W || G ]] is an epimorphism, and $V$ is additive, then the
induced map \mapx[[ \Mor_{\smallcat{$k$-functors}}(V,W) ||
                    \Mor_{\smallcat{$k$-functors}}(V,G) ]] is surjective.
\item\label{additive-any}
     If $V$ is additive, then $\Ext^1_s(V,F) = \Ext^1(V,F)$ for all $F$.
\item\label{Z-sheaf} If $F$ is a sheaf then $\Ext^1(\Z,F) = 0$.
\item\label{Gm-Ga} If we have an exact sequence
\ses{U_1}{U_2}{P%
}in which $U_1$ and $U_2$ are additive, then $\Ext^1(\Gm, P) = 0$.
\item\label{nilpotent-discrete} $\Ext^1(I,D) = 0$ if $I$ is subnilpotent and
$D$ is discrete.
\item\label{additive-additive-char0} $\Ext^1(U,V) = 0$ if $U$ and $V$ are
additive and $k$ has characteristic zero.
\end{alphalist}
\end{prop}

\begin{proof}
{\bf (a):\ }
First we prove this when $V = \Ga^\alpha$ for some $\alpha$.  If
$\alpha < \infty$, the claim is immediate.  Otherwise, the essential point is
that in a commutative diagram
we can (exercise) fill in a dotted arrow as shown.

Now suppose that $V$ is arbitrary.  From what we have just shown, it follows
that $\Ext^1_s(\Ga^\alpha,F) = \Ext^1(\Ga^\alpha,F)$ for all $F$.  In turn,
this implies that $V \iso \Ga^{\kern1pt\beta}$ in \cat{$k$-functors}, for some
$\beta$.
Hence (a) holds when $V$ is arbitrary.

\vspace{0.1in}

\par\noindent{\bf (b):\ } follows immediately from (a).

\vspace{0.1in}

\par\noindent{\bf (c):\ }
We have to show that if \mp[[ \pi || H || \Z ]] is an epimorphism in
\cabkf, and $H$ is a sheaf, then
$\pi$ splits.  For each $n \in \Z$, let $y_n \in \Z(k)$ correspond to the
constant map \mapx[[ \Spec(k) || \Z ]] of topological spaces with value $n$,
and choose $x_1 \in H(k)$ such that $\pi(x_1) = y_1$.  For each $n \in \Z$,
define $x_n \in H(k)$ to be $n x_1$.  Define \mp[[ \sigma || \Z || H ]] as
follows.  For any ring $B$, an element $\lambda \in \Z(B)$ corresponds to a
locally constant map \mapx[[ \Spec(B) || \Z ]] of topological spaces, and
therefore we may write
$B = B_1 \manytimes B_n$ in such a way that $\lambda$ is induced by
$(y_{r_1},\ldots,y_{r_n})$ for suitable $\vec r1n \in \Z$.  Since $H$ is a
sheaf for the Zariski topology, there is a unique element
$x_{r_1,\ldots,r_n} \in H(k^n)$ whose image in $H(k)$ under the \th{i}
projection map is $x_{r_i}$.  Now set $\sigma(\lambda)$ equal to the image
of $x_{r_1,\ldots,r_n}$ under the canonical map
\mapx[[ H(k^n) || H(B_1 \manytimes B_n) ]].  This defines $\sigma$, and thus
proves that $\pi$ splits.

\vspace{0.1in}

\par\noindent{\bf (d):\ }
First suppose that $\CHAR(k) = 0$.  Then $P \iso \Ga^\alpha$ for some $\alpha$.
If $\alpha < \infty$, the statement follows from
[\Lcitemark 10\Rcitemark \ 5.1.1(i)].
We have $\Ext^1(\Gm,\Ga^\infty) = \Ext^1_s(\Gm,\Ga^\infty)$.
Therefore an extension of $\Gm$ by $\Ga^\infty$ corresponds to a symmetric
factor system \mp[[ f || \Gm \times \Gm || \Ga^\infty ]].  For $n \gg 0$, we
can find a morphism \mp[[ f_n || \Gm \times \Gm || \Ga^n ]] through which $f$
factors.  But then $f_n$ is a symmetric factor system, and so $f_n$ is
trivial, since we already know that $\Ext^1(\Gm,\Ga^n) = 0$.  Hence $f$ is
trivial.  Hence $\Ext^1(\Gm,\Ga^\infty) = 0$.

Now suppose that $k$ has characteristic $p > 0$.
Since $\Gm$ and $P$ are sheaves, it suffices to show that
$\Extfpqc^1(\Gm,P) = 0$.  For $n$ sufficiently large, multiplication by
$p^n$ gives a zero map from $P$ to $P$.  It follows from the fpqc-exact
sequence
\sesmapsone{\mu_{p^n}}{}{\Gm}{{p^n}}{\Gm%
}that it is enough to show $\Hom(\mu_{p^n}, P) = 0$.  Let
\mp[[ f || \mu_{p^n} || P ]] be a morphism.  Let $H$ be the fiber product of
$\mu_{p^n}$ and $U_2$ over $P$.  Then we have an exact sequence:
\diagramno{(*)}{\rowfive{0}{U_1}{H}{\mu_{p^n}}{1.}%
}We will show that this sequence splits.  We can do this by showing that
$\Ext^1(\mu_{p^n},U_1) = 0$, but by the definition of additive, it is clearly
enough to show that $\Ext^1(\mu_{p^n}, \Ga^\alpha) = 0$ for all $\alpha$.
Arguing as in the characteristic zero case, one sees further that it is further
enough to show that $\Ext^1(\mu_{p^n}, \Ga) = 0$.  This is a special case of
\Lcitemark 10\Rcitemark \Rspace{}\ 5.1.1(d).  Hence $(*)$ splits.
Hence there exists a morphism \mp[[ \sigma || \mu_{p^n} || U_2 ]] such that
$\pi \circ \sigma = f$, where \mp[[ \pi || U_2 || P ]] is the given map.
Now I claim that $\sigma = 0$.  For this (arguing as above), it is enough
to show that $\Hom(\mu_{p^n}, \Ga) = 0$.  It is enough to do this when
$k = k^a$, and then the statement is well-known.  Hence $f = 0$.  Hence
$\Hom(\mu_{p^n}, P) = 0$, which completes the proof.

\vspace{0.1in}

\par\noindent{\bf (e):\ }
Let
\ses{D}{L}{I%
}be an exact sequence of (abelian group)-valued $k$-functors.
Define a $k$-functor $I'$ by
$I'(B) = \Ker[L(B)\ \rightarrow\ L(\RED{B})]$.  Then $I'$ defines a
splitting of the sequence.

\vspace{0.1in}

\par\noindent{\bf (f):\ }  It suffices to show that
$\Ext^1(\Ga^\alpha, \Ga^{\kern1pt\beta}) = 0$ for all $\alpha, \beta$.
Moreover, since
$\Ext^1$ converts a coproduct in the first variable into a product,
\WMAT\ $\alpha = 1$.  By (b), it suffices to show that
$\Ext^1_s(\Ga, \Ga^{\kern1pt\beta}) = 0$.  Arguing as in (d), it suffices to
show
that $\Ext^1_s(\Ga,\Ga^n) = 0$, and moreover we may as well take $n=1$.
Suppose we have an exact sequence
\sesdot{\Ga}{X}{\Ga%
}By (\Lcitemark 10\Rcitemark \ 3.9 ter.), $X \iso \Ga^2$.  But
(in characteristic zero) morphisms from $\Ga^n$ to $\Ga^m$ are in bijective
correspondence with vector space homomorphisms from $k^n$ to $k^m$, so the
sequence splits.  \qed
\end{proof}

\block{Functorial structure of units in a tensor product}

The main purpose of this section is to prove \pref{tori-result-generalized},
which generalizes (\Lcitemark 8\Rcitemark \ 4.5).  The
preparatory lemmas are similar to those in
(\Lcitemark 8\Rcitemark \ \S4),
and we shall omit their proofs if the proofs of the corresponding statements
in\Lspace \Lcitemark 8\Rcitemark \Rspace{} carry over with minor changes.

\begin{lemma}\label{unit-lemma-generalized}
Let $k$ be an algebraically closed field.  Let $A$ and $B$ be reduced rings
containing $k$, having connected spectra.  Let $u \in A \o*_k B$ be a unit.
Then $u = a \o* b$ for some units $a \in A$ and $b \in B$.
\end{lemma}

\begin{proof}
The statement generalizes (\Lcitemark 8\Rcitemark \ 4.2), but we
give a new and simpler proof, due to Guralnick.

Let $X$ be the set of maximal ideals of $A$, and let $Y$ be the set of maximal
ideals of $B$.  Let $x_0 \in X$, $y_0 \in Y$.  We will prove the lemma by
showing that for all $x \in X$, $y \in Y$, we have:
$$u(x,y)\ =\ {u(x,y_0) u(x_0,y) \over u(x_0, y_0)}.\eqno(*)$%
$For this we may suppose that $k$ is uncountable.
By a variant of a result of Roquette (see
\Lcitemark 7\Rcitemark \Rspace{}\ 1.5) the group $B^*/k^*$ is \FG, so
$F = \setof{f \in B^*: f(y_0) = 1}$ is countable.  For each $f \in F$, let
$$Q(f)\ =\ \setof{x \in X: u(x,y) = u(x,y_0)f(y)\hbox{\ for all\ } y \in Y}.$%
$Then $Q(f)$ is a closed subset of $X$.  For any given $x \in X$, the
function on $Y$ given by $y \mapsto u(x,y)/u(x,y_0)$ sends $y_0$ to $1$
and so lies in $F$.  Hence $X = \cup_{f \in F}Q(f)$.

Since $k$ is uncountable, it follows that if $I$ is an irreducible component
of $X$, then $I \IN Q(f)$ for some $f \in F$.  Hence for any fixed $y \in Y$,
the function \mp[[ g_y || X || k ]] given by $x \mapsto u(x,y)/u(x,y_0)$ is
constant on each irreducible component of $X$.  Since $X$ is connected, $g_y$
is constant.  Then $u(x,y) = u(x,y_0)g_y(x_0)$, which proves $(*)$.  \qed
\end{proof}

\begin{corollary}\label{second-generalized}
Let $k$ be an algebraically closed field.  Let $A$ and $B$ be rings containing
$k$, having connected spectra.  Assume that $A$ is reduced.  Let $\lfM \IN A$
be a maximal ideal such that $A/\lfM = k$.  Then $(A \o*_k B)^*$ is the direct
sum of the two subgroups $A^*B^*$ and $1 + (\lfM \o* \Nil(B))$%
.\footnote{Statements 4.3 and 4.4 from\Lspace \Lcitemark 8\Rcitemark \Rspace{}
should
also have the hypothesis that $A/\lfM = k$.}
\end{corollary}

\begin{corollary}\label{third-generalized}
Let $k$ be an algebraically closed field.  Let $A$ and $B$ be rings
containing $k$.  Assume that $A$ is reduced and has a connected spectrum.
Let $\lfM \IN A$ be a maximal ideal such that $A/\lfM = k$.
\begin{itemize}
\item For any decomposition
$B = B_1 \manytimes B_n$, there is a subgroup:
$$\mu(\vec B1n) = \oplus_{i=1}^n [A^*B_i^* \o+ (1 + (\lfM \o* \Nil(B_i))]$%
$of $(A \o*_k B)^*$.
\item For any $x \in (A \o*_k B)^*$, there exists a decomposition
$B = B_1 \manytimes B_n$ such that $x \in \mu(\vec B1n)$.
\item If $B$ has only finitely many idempotent elements (e.g.\ if $B$ is
noetherian), we can write $B = B_1 \manytimes B_n$ for
rings $B_i$ having connected spectra.  Then $\mu(\vec B1n) = (A \o*_k B)^*$.
\end{itemize}
\end{corollary}

Let $k$ be a field, and let $S$ be a $k$-scheme of finite type.  Let
$F = \HOM(S,\Gm)$.  Then
$$F(B)\ =\ \Gamma(S \times \Spec(B), \O^*_{S \times \Spec(B)})
      \ =\ [\Gamma(S,\O_S) \o*_k B]^*,$%
$by (\Lcitemark 6\Rcitemark \ 9.3.13 (i)).
In particular, if $S = \Spec(A)$, then $F(B) = (A \o*_k B)^*$.  The next
theorem gives an abstract description of $F$, and thus (in effect) a
description of how units in a tensor product $A \o*_k B$ vary as $B$ varies.
First, for convenience, we encapsulate the following definition:

\begin{definition}
Let $k$ be a field.  A $k$-scheme $S$ is {\it geometrically stable\/} if $(1)$
it is of finite type, and $(2)$ every irreducible component of $\RED{S}$ is
geometrically integral and has a rational point.
\end{definition}

\begin{theorem}\label{tori-result-generalized}
Let $k$ be a field.  Define an (abelian group)-valued $k$-functor $F$ to be of
type $(*)$ if there exist exact sequences
\seslabcomma{R}{F}{I}{\dag%
}\seslab{\Gm^r \times U}{R}{L}{\dag\dag%
}in \cabkf, in which
$r \geq 0$, $U$ is pseudoadditive, $I$ is subnilpotent, and $L$ is discrete
and \FG.
\par\noindent{\bf\rm (a):\ }
Let $S$ be a geometrically stable $k$-scheme.  Then $\HOM(S,\Gm)$ is of type
$(*)$ and we have $r =$ the number of connected components of $S$.  Also $U$
is additive of dimension $\dim_k \Nil[\Gamma(S,\O_S)]$.  Also,
$I$ is nilpotent, $L$ is free, and $(\dag\dag)$ splits.
\par\noindent{\bf\rm (b):\ }
Let $S$ and $T$ be geometrically stable $k$-schemes, and let
\mp[[ f || S || T ]] be a dominant morphism of $k$-schemes.
Then the cokernel of $\HOM(f,\Gm)$ is of type $(*)$.  Moreover, $r$ equals
the number of connected components of $S$ minus the number of connected
components of $T$.
\end{theorem}

\begin{corollary}\label{tori-quotient-sheaf}
Let $S$ and $T$ be geometrically stable $k$-schemes, and let
\mp[[ f || S || T ]] be a dominant morphism of $k$-schemes.  Let $Q$ be the
cokernel of $\HOM(f,\Gm)$.  Then the canonical map \mapx[[ Q || Q^+ ]] is a
monomorphism, and $Q|_{\smallcat{reduced $k$-algebras}}$ is a sheaf, in the
sense that if \mp[[ p || B || C ]] is a faithfully flat homomorphism of
reduced $k$-algebras, then $\psi_{Q,p}$ (see p.\ \pageref{Psi-place}) is
bijective.
\end{corollary}

\begin{remarks}
\
\begin{romanlist}
\item In part (a) of the theorem, one can choose $(\dag)$ so that it splits
      if $S$ is reduced, but probably not in general.
\item In part (b), the sequence $(\dag\dag)$ does not always split.  For
      an example, take $k$ to be an imperfect field of characteristic $p$,
      let $u \in k - k^p$, and let $f$ be $\Spec$ of the ring map
      \mapx[[ k[t,t^{-1}] || k[x,x^{-1}] \times k ]], given by
      $t \mapsto (x^p,u)$.
\item The hypothesis that the schemes in the theorem be geometrically stable
      can be weakened slightly, as is indicated in the proof.  They presumably
      can be weakened further, but we do not know what is possible in this
      direction.
\item\label{tori-result-remark-three}
      We suspect that $I$ in part (b) of the theorem is a sheaf (and thus
      satisfies the definition of {\it nilpotent}).  If true, this would imply
      (in the corollary) that the cokernel of $\HOM(f,\Gm)$ is a sheaf.  To
      prove that $I$ is a sheaf, it would be sufficient (at least in the
      case where $S$ and $T$ are connected) to show that if $C$ is
      a subalgebra of a reduced $k$-algebra $A$, then the
      (abelian group)-valued $k$-functor given by
      $$B\ \mapsto\ {1 + A \o* \Nil(B) \over 1 + C \o* \Nil(B)}$%
$      is a sheaf.
\item In part (b), we have $U$ pseudoadditive with $n=2$, as in the definition
      of pseudoadditive.  However, it is conceivable that $U$ is always
      additive.
\end{romanlist}
\end{remarks}

\begin{proofnodot}
(of \ref{tori-result-generalized}.)
The hypothesis that $S$ be geometrically stable is chosen for simplicity
and we note here some consequences which are in fact sufficient to prove the
theorem:
\vspace{0.05in}
\par\circno{A}:\ every connected component of $S$ has a rational
                          point;
\vspace{0.05in}
\par\circno{B}:\ $\RED{S}$ is geometrically reduced
                          [by\Lspace \Lcitemark 5\Rcitemark \Rspace{}\
4.6.1(e)].
\vspace{0.05in}
\par\noindent Moreover, \circno{A}\ and \circno{B}\ also imply:
\vspace{0.05in}
\par\circno{C}:\ every connected component of $S$ is geometrically
                          connected [by\Lspace \Lcitemark 5\Rcitemark
\Rspace{}\ 4.5.14];
\vspace{0.05in}
\par\circno{D}:\ if $Q$ is a connected component of $S$, then $k$ is
                          integrally closed in $\Gamma(\RED{Q})$.
\par\noindent All of these comments apply equally to $T$.

Now we want to reduce to the affine case.  This is not literally possible,
because $\Gamma(S)$ need not be \FG\ as a $k$-algebra.  What we can do is
reformulate the theorem in terms of a certain class of $k$-algebras.  This
class is chosen simply to serve the needs of the proof: a $k$-algebra
is {\it good\/} if it is of the form $\Gamma(S)/N$, where $S$ is a
geometrically stable $k$-scheme and $N \IN \Gamma(S)$ is a nilpotent ideal.
Here is a reformulation of the theorem in terms of good $k$-algebras:
\begin{quote}
\par\noindent{\bf\rm (a):\ }
Let $A$ be a good $k$-algebra.  Then $B \mapsto (A \o* B)^*$ is of type
$(*)$ and we have $r =$ the number of connected components of $\Spec(A)$,
$\dim(U) = \dim_k \Nil(A)$.  Also, $L$ is free and $I$ is nilpotent.
\par\noindent{\bf\rm (b):\ }
Let \mp[[ \phi || C || A ]] be a homomorphism of good $k$-algebras.
Assume that $\Ker(\phi)$ is nilpotent.  Then the cokernel of the morphism from
$B \mapsto (C \o* B)^*$ to $B \mapsto (A \o* B)^*$ is of type $(*)$.
Moreover, $r$ equals the number of connected components of $\Spec(A)$ minus
the number of connected components of $\Spec(C)$.
\end{quote}

Since the map \mapx[[ (C \o* B)^* || (C/\Ker(\phi) \o* B)^* ]] is surjective,
we may reduce to the case where $\phi$ is {\it injective}.  It was the need for
this reduction which lead to the introduction of good $k$-algebras in the
proof.

We proceed to build a diagram involving (abelian group)-valued $k$-functors,
which we associate to $A$, and which is functorial in $A$.

Let $G_A$ be given by $G_A(B) = (\RED{A} \o* B)^*$.  Write
$\RED{A} = A_1 \manytimes A_r$, where $\vec A1r$ have connected spectra.
We can identify $A \o* B$ with $(A_1 \o* B) \manytimes (A_r \o* B)$.  Let
$F_A$ be the sheaf associated to the subfunctor of $G_A$ given by
$B \mapsto
   \setof{(a_1 \o* b_1, \ldots, a_r \o* b_r): a_i \in A_i^*, b_i \in B^*}$.
Let $E_A$ be the subfunctor of $F_A$ given by
$E_A(B) = \setof{(\vec b1r): \vec b1r \in B^*}$.  Let
$D_A = F_A/E_A$.  Let $I_A = G_A/F_A$.

Define $H_A$ by $H_A(B) = (A \o* B)^*$.
Let \mp[[ p || H_A || G_A ]] be the canonical map, and let $U_A$ be its kernel.
We have $U_A(B) = 1 + \Nil(A) \o* B$.  For each $n \in \N$, let
$U_A^n$ be given by $U_A^n(B) = 1 + \Nil(A)^n \o* B$.

Here is the diagram of (abelian group)-valued $k$-functors which we have
built:

\diagramx{&&&&0\cr
          &&&&\mapS{}\cr
          && 0 &&U_A & = \kern10pt
                         \hbox to 0pt{$U_A^1 \kern10pt \NI \kern10pt U_A^2
                                      \kern10pt \NI \kern10pt \cdots$}\cr
          && \mapS{} && \mapS{}\cr
          && E_A && H_A \cr
          && \mapS{} && \mapS{}\cr
          \rowfive{0}{F_A}{G_A}{I_A}{0}\cr
          && \mapS{} && \mapS{} && \vbox to 0pt{\box5}\cr
          && D_A && 0\cr
          && \mapS{}\cr
          && 0}

We proceed to analyze the various components of this diagram.  In particular,
we will show that \circno1\ $G_A \iso F_A \times I_A$ and $I_A$ is nilpotent,
\circno2\ $E_A \iso \Gm^r$,
where $r$ is the number of connected components of $\Spec(A)$, \circno3\ $D_A$
is represented by a constant group scheme (corresponding to a \FG\ free abelian
group), and that \circno4\ $U_A$ is additive.  For these purposes,
\WMAT\ $\Spec(A)$ is connected.  Then $E_A(B) = B^*$.

\vspace{0.05in}
\par\noindent\circno1\ Let
$\lfM_A \IN \RED{A}$ be a maximal ideal such that $A/\lfM_A = k$.  (This is
possible by \circno{A}.)  Define a
subfunctor $I'_A$ of $G_A$ by $I'_A(B) = 1 + \lfM_A \o* \Nil(B)$.
Let \mp[[ \psi || F_A \o+ I'_A || G_A ]] be the canonical map.  We will show
that $\psi$ is an isomorphism.  If $k$ is algebraically closed, this follows
from \pref{third-generalized}.  But both the source and the target of $\psi$
are sheaves, so it follows (using \circno{B}\ and \circno{C}) that $\psi$ is an
isomorphism for any $k$.  Thus $I'_A \iso I_A$, so $I_A$ is nilpotent.

\vspace{0.05in}
\par\noindent\circno2\ We have $E_A \iso \Gm$.

\vspace{0.05in}
\par\noindent\circno3\ Write $A = \Gamma(S)/N$, as in the definition of
good.  Let $\shN$ be the nilradical of $S$.  Then we have an exact sequence:
\les{H^0(S,\shN)}{H^0(S,\O_S)}{H^0(\RED{S},\O_{\RED{S}})%
}and so $\RED{A}$ is a subring of $H^0(\RED{S},\O_{\RED{S}})$.  Since $k$ is
integrally closed in $H^0(\RED{S},\O_{\RED{S}})$ by \circno{D}, it follows by
(\Lcitemark 7\Rcitemark \ 1.5),that $D_A(k)$ is free abelian of
finite rank.  In fact, $D_A$ is represented by the corresponding constant group
scheme.

\vspace{0.05in}
\par\noindent\circno4\ For each $n$ we have
$${U_A^n \over U_A^{n+1}}(B)
  \ =\ {1 + \Nil(A)^n \o* B \over 1 + \Nil(A)^{n+1} \o* B}
  \ \iso\ {\Nil(A)^n \over \Nil(A)^{n+1}} \o* B$%
$as (abelian group)-valued $k$-functors.
Hence $U_A$ is additive.

Now we describe $H_A$, making a number of non-canonical choices.
We have $F_A \iso E_A \times D_A$ non-canonically, e.g.\ by
(\ref{ext-results}\ref{Z-sheaf}), but it is easily proved directly.  Hence
$F_A \iso \Gm^r \times \Z^n$ for some $n$.  Also, we have shown that
$G_A \iso F_A \times I_A$.  Therefore we have an exact sequence
\sesmapsdot{U_A}{}{H_A}{q}{\Gm^r \times \Z^n \times I_A%
}Let $M = q^{-1}(\Gm^r \times \Z^n)$.  Then we have an exact sequence
\sesdot{U_A}{M}{\Gm^r \times \Z^n%
}By (\ref{ext-results}\ref{Z-sheaf}\ref{Gm-Ga}) this sequence splits.
This proves (a).

Now, to prepare for proving (b),  we analyze the functorial behavior of each
basic component of the big diagram shown above.  Let \mp[[ \phi || C || A ]] be
an injective homomorphism of good $k$-algebras.

First we analyze $E_\phi$.  It is a monomorphism,
corresponding to a map \mapx[[ \Gm^{r_1} || \Gm^{r_2} ]], for some $r_1$ and
$r_2$, which is given by an $r_2 \times r_1$ matrix of $0$'s and $1$'s.
The cokernel of $E_\phi$ is isomorphic to $\Gm^{r_2 - r_1}$.

Now we analyze $D_\phi$.  Let us show that $D_\phi$ is a monomorphism.
Since its source and target are constant sheaves, it suffices to show that
$D_\phi(k^a)$ is injective.  The assertion then boils down to showing that
if one has a dominant morphism \mp[[ \psi || V || W ]] of reduced schemes of
finite type over an algebraically closed field $k$, and $W$ is connected, and
\mp[[ g || W(k) || k ]] is a non-constant regular function, then
$g \circ \psi(k)$ is not constant on each connected component of $V$.  This is
clear, so $D_\phi$ is a monomorphism.  The cokernel of $D_\phi$ is the
constant sheaf associated to a \FG\ abelian group.

We show that $I_\phi$ is a monomorphism.  In the process of doing so, we
justify remark \pref{tori-result-remark-three} from
p.\ \pageref{tori-result-remark-three}.  Also, once we know that $I_\phi$ is a
monomorphism, it will follow immediately that $\Coker(I_\phi)$ is subnilpotent.
We may assume that $\Spec(C)$ is connected.

We have a canonical map \mapx[[ 1 + \RED{C} \o* \Nil(B) || I_C(B) ]], and
likewise for $A$.  From our discussion of $I'$, it is clear that these maps
are surjective.  Letting $X_C(B)$ and $X_A(B)$ denote their kernels, we have a
commutative diagram with exact rows:
\diagramx{\sesonerow{X_C(B)}{1 + \RED{C} \o* \Nil(B)}{I_C(B)}\cr
          && \mapS{} && \mapS{} && \mapS{}\cr
          \sesonerowdot{X_A(B)}{1 + \RED{A} \o* \Nil(B)}{I_A(B)}}

We describe $X_C(B)$.
Let $x \in X_C(B)$.  Locally (for the fpqc topology) on $B$, we may write
$x = c \o* b = 1 + \sum_{i=1}^s c_i \o* b_i$, where $c \in \RED{C}^*$,
$b \in B^*$, $c_i \in \RED{C}$, and $b_i \in \Nil(B)$.  It follows that $b$
must lie in the $k$-linear span of $1$ and the $b_i$, and in particular that
we may write $b = 1 + n$, where $n \in \Nil(B)$.  Passing to
$\RED{C} \o* \RED{B}$, we see then that $c = 1$.  Hence $x = 1 + n$.  From
this it follows that $X_C(B) = 1 + \Nil(B)$.

Similarly, we have
$X_A(B) = (1+\Nil(B)) \manytimes (1+\Nil(B))$, with one copy for each
connected component of $\Spec(A)$.  It follows (details omitted) that the
canonical map
\dmapx[[ {X_A(B) \over X_C(B)} ||
         {1 + \RED{A} \o* \Nil(B) \over 1 + \RED{C} \o* \Nil(B)} ]]%
is injective, and hence that $I_\phi$ is a monomorphism.

Assume now that $\Spec(A)$ and $\Spec(C)$ are connected.
Let $\lfM_C$ be the preimage of $\lfM_A$ under the map
\mapx[[ \RED{C} || \RED{A} ]] induced by $\phi$.  From our discussion of $I'$,
it is clear that $\Coker(I_\phi)$ is isomorphic to the cokernel of the morphism
given at $B$ by \mapx[[ 1 + \lfM_C \o* \Nil(B) || 1 + \lfM_A \o* \Nil(B) ]].
In turn this implies remark (iii).

We have exact sequences
\sescomma{\Coker(U_\phi)}{\Coker(H_\phi)}{\Coker(G_\phi)%
}\sescomma{\Coker(F_\phi)}{\Coker(G_\phi)}{\Coker(I_\phi)%
}\sesdot{\Coker(E_\phi)}{\Coker(F_\phi)}{\Coker(D_\phi)%
}Since $\Coker(U_\phi)$ is clearly pseudoadditive, part (b) of the theorem
follows from these sequences and (\ref{ext-results}\ref{Gm-Ga}).  \qed
\end{proofnodot}

\block{Line bundles becoming trivial on pullback by a nilimmersion}

If $X$ is a $k$-scheme, we let $\PIC(X)$ denote the $k$-functor given by
$B \mapsto \Pic(X \times_k \Spec(B))/\Pic(B)$.  Then $\PIC$ itself defines a
functor whose source is \opcat{$k$-schemes}.  If \mp[[ f || X || Y ]] is a
morphism of $k$-schemes of finite type, such that $X$ and $Y$ each have a
rational point, then $\Ker[\PIC(f)]$ is isomorphic to the $k$-functor given by
$$B \mapsto \Ker[\Pic(Y \times_k \Spec(B))
  \ \rightarrow\ \Pic(X \times_k \Spec(B))].$$

\begin{theorem}\label{kernel-pic-nilimmersion}
Let $k$ be a field, and let $X$ be a geometrically stable $k$-scheme.
Let \mp[[ i || X_0 || X ]] be a nilimmersion, such that the ideal sheaf $\shN$
of $X_0$ in $X$ has square zero.  Then there is an exact sequence
of (abelian group)-valued $k$-functors
\sescomma{D \o+ I}{P}{\Ker[\PIC(i)]%
}in which $D$ is discrete and \FG, $I$ is subnilpotent, and $P$ is
pseudoadditive.
\end{theorem}

\begin{remarks}
\
\begin{alphalist}
\item If $X$ is affine, $\Ker[\PIC(i)] = 0$.
\item If $X$ is proper over $k$, $D = 0$ and $I = 0$, so
      $\Ker[\PIC(i)] \iso P$.  Also, $P$ is additive and finite-dimensional.
      One way to get examples is to take $k$
      to be algebraically closed, $Y$ to be a projective variety over $k$, and
      $\shM$ to be a coherent $\O_Y$-module with $H^1(Y,\shM) \not= 0$.  Make
      $\O_Y \o+ \shM$ into a coherent $\O_Y$-algebra via the rule $\shM^2 = 0$.
      Let $X = \SPEC(\O_Y \o+ \shM)$, and let \mp[[ i || \RED{X} || X ]] be the
      inclusion.  Then $\Ker[\PIC(i)] \iso \Ga^\alpha$, where
      $\alpha = h^1(Y,\shM)$.
\item In the non-affine, non-proper case, we have not determined exactly what
      can happen.  In particular, we do not know if $D$ can be nonzero.
      If $k$ has characteristic zero, then $D \iso \Z^n$ for some $n$.  If $k$
      has characteristic $p > 0$, then $D \iso \Z^n \o+ (\Z/p\Z)^m$ for some
      $n$ and some $m$.
\item Conceivably the theorem holds without the assumption that $\shN^2 = 0$.
      To prove this, one would at least have to understand
      $\Coker[\PIC(i)]$ in the case where $\shN^2 = 0$, which we do not.
\item We will find an exact sequence
\Rowsix{0}{U_1}{U_2}{\Ga^{\kern1pt\beta}}{P}{0%
}in which $U_1$ and $U_2$ are additive, for some $\beta$.  This is
stronger than saying that $P$ is pseudoadditive.  Perhaps $P$ is always
additive.
\end{alphalist}
\end{remarks}

\begin{proof}
We have an exact sequence of sheaves of abelian groups on $X$
\Rowfive{0}{\shN}{\O_X^*}{\O_{X_0}^*}{1,%
}and thus an exact sequence of abelian groups
\splitdiagram{H^0(X,\O_X^*)&\mapE{}&H^0(X_0,\O_{X_0}^*)%
}{\mapE{}&H^1(X,\shN)&\mapE{}&
\Ker[\Pic(X)\ \rightarrow\ \Pic(X_0)]&\mapE{}&0.%
}In fact, everything is functorial in $k$, and we thus obtain an exact
sequence of (abelian group)-valued $k$-functors
\splitdiagram{\HOM(X,\Gm)&\mapE{}&\HOM(X_0,\Gm)&\mapE{}&
\Ga^{\kern1pt\beta}}{\mapE{}&
\Ker[\PIC(X)\ \rightarrow\ \PIC(X_0)]&\mapE{}&0,%
}where $\beta = \dim_k[H^1(X,\shN)]$.
Let $K = \Ker[\PIC(X)\ \rightarrow\ \PIC(X_0)]$,
$L = \Coker[\HOM(X,\Gm)\ \rightarrow\ \HOM(X_0,\Gm)]$.

We have an exact sequence
\sesdot{L}{\Ga^{\kern1pt\beta}}{K%
}According to (\ref{tori-result-generalized}b), there are exact sequences
\sescomma{P'}{R}{Q%
}\sescomma{R}{L}{I%
}in which $Q$ is discrete and \FG, $P'$ is additive, and $I$ is subnilpotent.
By definition,
$\Coker[ P'\ \rightarrow\ \Ga^{\kern1pt\beta}]$ is pseudoadditive.  Hence by
(\ref{ext-results}\ref{nilpotent-discrete}), we have an exact sequence
\sescomma{Q \times I}{P}{K%
}in which $P$ is pseudoadditive.  \qed
\end{proof}

\begin{corollary}\label{sheaf-kernel-pic-nilimmersion}
Let $k$ be a field, and let $X$ be a geometrically stable $k$-scheme.
Let \mp[[ i || X_0 || X ]] be a nilimmersion.  Let
$F = \Ker[\PIC(i)]$.  Let \mp[[ p || B || C ]] be a faithfully flat ring
homomorphism.
\begin{alphalist}
\item If $B$ and $C$ are reduced, then the canonical map
\mapx[[ F(B) || F(C) ]] is injective.
\item Assume that $k$ has characteristic zero and that the ideal sheaf of $X_0$
in $X$ has square zero.  Assume that $B$ is normal and
that $C$ is \'etale over $B$.  Then $\Psi_{F,p}$ (see p.\ \pageref{Psi-place})
is bijective.
\end{alphalist}
\end{corollary}

\begin{proof}
Let $\shN$ be the ideal sheaf of $X_0$ in $X$.  First suppose that
$\shN^2 = 0$.  Then (a) follows from \pref{kernel-pic-nilimmersion}.  For (b),
let $S = \Spec(B)$.  Since $\CHAR(k) = 0$, $D$ is torsion-free.  Therefore
it suffices to show that in the category of
(abelian group)-valued $k$-functors
which are sheaves for the \'etale topology, the sequence
\ses{\Z^n}{\Ga^\alpha}{(\Ga^\alpha/\Z^n)^+%
}is exact when evaluated at $B$.  (Here we let $(\Ga^\alpha/\Z^n)^+$ denote the
quotient in this category.)  Since $H^1(\et{S}, \Z) = 0$ by
(\Lcitemark 1\Rcitemark \ 3.6(ii)), we are done.

Now we prove the general case of (a).  For each $m$, let $X_m$ be the closed
subscheme of $X$ defined by $\shN^m$.  Choose $n \in \N$ so that $\shN^n = 0$.
Let $K_m = \Ker[\PIC(X_{m+1})\ \rightarrow\ \PIC(X_m)]$.
Let $F_m = \Ker[\PIC(X_m)\ \mapE{}\ \PIC(X_0)]$.  We have an exact sequence
\Rowfour{0}{K_m}{F_{m+1}}{F_m.%
}Then \mapx[[ K_m(B) || K_m(C) ]] is injective by
\pref{kernel-pic-nilimmersion}.
By induction on $m$, it follows that \mapx[[ F_m(B) || F_m(C) ]] is injective
for all $m$.  Taking $m = n$, we get (a).  \qed
\end{proof}

\begin{problem}
Is the functor $F$ a sheaf?
\end{problem}
\section*{References}
\addcontentsline{toc}{section}{References}
\ \par\noindent\vspace*{-0.25in}
\hfuzz 5pt

\bgroup\Resetstrings%
\def\Loccittest{}\def\Abbtest{ }\def\Capssmallcapstest{}\def\Edabbtest{
}\def\Edcapsmallcapstest{}\def\Underlinetest{ }%
\def\NoArev{1}\def\NoErev{0}\def\Acnt{1}\def\Ecnt{0}\def\acnt{0}\def\ecnt{0}%
\def\Ftest{ }\def\Fstr{1}%
\def\Atest{ }\def\Astr{Artin\Revcomma M\Initper }%
\def\Ttest{ }\def\Tstr{Faisceaux constructibles, cohomologie d'un courbe
alg\`ebrique, \EXPOSE\ IX in {\it\SGA} (SGA 4)}%
\def\Stest{ }\def\Sstr{Lecture Notes in Mathematics}%
\def\Itest{ }\def\Istr{Spring\-er-Ver\-lag}%
\def\Ctest{ }\def\Cstr{New York}%
\def\Vtest{ }\def\Vstr{305}%
\def\Dtest{ }\def\Dstr{1973}%
\def\Ptest{ }\def\Pcnt{ }\def\Pstr{1--42}%
\def\Qtest{ }\def\Qstr{access via "artin constructible sheaves"}%
\Refformat\egroup%

\bgroup\Resetstrings%
\def\Loccittest{}\def\Abbtest{ }\def\Capssmallcapstest{}\def\Edabbtest{
}\def\Edcapsmallcapstest{}\def\Underlinetest{ }%
\def\NoArev{1}\def\NoErev{0}\def\Acnt{1}\def\Ecnt{0}\def\acnt{0}\def\ecnt{0}%
\def\Ftest{ }\def\Fstr{2}%
\def\Atest{ }\def\Astr{Artin\Revcomma M\Initper }%
\def\Ttest{ }\def\Tstr{The implicit function theorem in algebraic geometry}%
\def\Btest{ }\def\Bstr{Algebraic Geometry (Bombay Colloquium, 1968)}%
\def\Itest{ }\def\Istr{Oxford Univ. Press}%
\def\Dtest{ }\def\Dstr{1969}%
\def\Ptest{ }\def\Pcnt{ }\def\Pstr{13--34}%
\def\Xtest{ }\def\Xstr{Theorem: $G/H$ is representable by an algebraic space.
repr-quot: section 7.}%
\def\Qtest{ }\def\Qstr{access via "artin implicit function theorem"}%
\def\Xtest{ }\def\Xstr{descent: p. 31}%
\Refformat\egroup%

\bgroup\Resetstrings%
\def\Loccittest{}\def\Abbtest{ }\def\Capssmallcapstest{}\def\Edabbtest{
}\def\Edcapsmallcapstest{}\def\Underlinetest{ }%
\def\NoArev{1}\def\NoErev{0}\def\Acnt{1}\def\Ecnt{0}\def\acnt{0}\def\ecnt{0}%
\def\Ftest{ }\def\Fstr{3}%
\def\Atest{ }\def\Astr{Bertin\Revcomma J\Initper \Initgap E\Initper }%
\def\Ttest{ }\def\Tstr{{\rm\tolerance=1000 Generalites sur les preschemas en
groupes,  \EXPOSE\ ${\rm VI}_{\rm B}$ in {\it\SGA} (SGA 3)}}%
\def\Stest{ }\def\Sstr{Lecture Notes in Mathematics}%
\def\Itest{ }\def\Istr{Spring\-er-Ver\-lag}%
\def\Ctest{ }\def\Cstr{New York}%
\def\Vtest{ }\def\Vstr{151}%
\def\Dtest{ }\def\Dstr{1970}%
\def\Ptest{ }\def\Pcnt{ }\def\Pstr{318--410}%
\def\Qtest{ }\def\Qstr{access via "bertin preschemas"}%
\def\Xtest{ }\def\Xstr{This includes the following result, announced by
Raynaud, but apparently not proved in this volume:  Theorem (11.11.1)  Let $S$
be a regular noetherian scheme of dimension $/leq 1$.  Let \map(\pi,G,S) be an
$S$-group scheme of finite-type.  Assume that $\pi$ is a flat, affine morphism.
 Then $G$ is isomorphic to a closed  subgroup of $\AUT(V)$, where $V$ is a
vector bundle on $S$.}%
\Refformat\egroup%

\bgroup\Resetstrings%
\def\Loccittest{}\def\Abbtest{ }\def\Capssmallcapstest{}\def\Edabbtest{
}\def\Edcapsmallcapstest{}\def\Underlinetest{ }%
\def\NoArev{1}\def\NoErev{0}\def\Acnt{1}\def\Ecnt{0}\def\acnt{0}\def\ecnt{0}%
\def\Ftest{ }\def\Fstr{4}%
\def\Atest{ }\def\Astr{Gabriel\Revcomma P\Initper }%
\def\Ttest{ }\def\Tstr{{\rm\tolerance=1000 Generalites sur les groupes
algebriques,  \EXPOSE\ ${\rm VI}_{\rm A}$ in {\it\SGA} (SGA 3)}}%
\def\Stest{ }\def\Sstr{Lecture Notes in Mathematics}%
\def\Itest{ }\def\Istr{Spring\-er-Ver\-lag}%
\def\Ctest{ }\def\Cstr{New York}%
\def\Vtest{ }\def\Vstr{151}%
\def\Dtest{ }\def\Dstr{1970}%
\def\Ptest{ }\def\Pcnt{ }\def\Pstr{287--317}%
\def\Qtest{ }\def\Qstr{access via "gabriel generalites"}%
\Refformat\egroup%

\bgroup\Resetstrings%
\def\Loccittest{}\def\Abbtest{ }\def\Capssmallcapstest{}\def\Edabbtest{
}\def\Edcapsmallcapstest{}\def\Underlinetest{ }%
\def\NoArev{1}\def\NoErev{0}\def\Acnt{2}\def\Ecnt{0}\def\acnt{0}\def\ecnt{0}%
\def\Ftest{ }\def\Fstr{5}%
\def\Atest{ }\def\Astr{Grothendieck\Revcomma A\Initper %
  \Aand J\Initper \Initgap A\Initper  Dieudonn\'e}%
\def\Ttest{ }\def\Tstr{El\'ements de \geom\ alg\'e\-brique IV (part two)}%
\def\Jtest{ }\def\Jstr{Inst. Hautes \'Etudes Sci. Publ. Math.}%
\def\Vtest{ }\def\Vstr{24}%
\def\Dtest{ }\def\Dstr{1965}%
\def\Qtest{ }\def\Qstr{access via "EGA4-2"}%
\Refformat\egroup%

\bgroup\Resetstrings%
\def\Loccittest{}\def\Abbtest{ }\def\Capssmallcapstest{}\def\Edabbtest{
}\def\Edcapsmallcapstest{}\def\Underlinetest{ }%
\def\NoArev{1}\def\NoErev{0}\def\Acnt{2}\def\Ecnt{0}\def\acnt{0}\def\ecnt{0}%
\def\Ftest{ }\def\Fstr{6}%
\def\Atest{ }\def\Astr{Grothendieck\Revcomma A\Initper %
  \Aand J\Initper \Initgap A\Initper  Dieudonn\'e}%
\def\Ttest{ }\def\Tstr{El\'ements de \geom\ alg\'e\-brique I}%
\def\Itest{ }\def\Istr{Springer-Verlag}%
\def\Ctest{ }\def\Cstr{New \snap York}%
\def\Dtest{ }\def\Dstr{1971}%
\def\Qtest{ }\def\Qstr{access via "EGA1"}%
\def\Astr{\Underlinemark}%
\Refformat\egroup%

\bgroup\Resetstrings%
\def\Loccittest{}\def\Abbtest{ }\def\Capssmallcapstest{}\def\Edabbtest{
}\def\Edcapsmallcapstest{}\def\Underlinetest{ }%
\def\NoArev{1}\def\NoErev{0}\def\Acnt{4}\def\Ecnt{0}\def\acnt{0}\def\ecnt{0}%
\def\Ftest{ }\def\Fstr{7}%
\def\Atest{ }\def\Astr{Guralnick\Revcomma R\Initper %
  \Acomma D\Initper \Initgap B\Initper  Jaffe%
  \Acomma W\Initper  Raskind%
  \Aandd R\Initper  Wiegand}%
\def\Ttest{ }\def\Tstr{The kernel of the map on Picard groups induced by a
faithfully flat homomorphism}%
\def\Rtest{ }\def\Rstr{preprint}%
\def\Qtest{ }\def\Qstr{access via "gang of four" or "guralnick jaffe raskind
wiegand"}%
\Refformat\egroup%

\bgroup\Resetstrings%
\def\Loccittest{}\def\Abbtest{ }\def\Capssmallcapstest{}\def\Edabbtest{
}\def\Edcapsmallcapstest{}\def\Underlinetest{ }%
\def\NoArev{1}\def\NoErev{0}\def\Acnt{1}\def\Ecnt{0}\def\acnt{0}\def\ecnt{0}%
\def\Ftest{ }\def\Fstr{8}%
\def\Atest{ }\def\Astr{Jaffe\Revcomma D\Initper \Initgap B\Initper }%
\def\Ttest{ }\def\Tstr{On sections of commutative group schemes}%
\def\Jtest{ }\def\Jstr{Compositio Math.}%
\def\Vtest{ }\def\Vstr{80}%
\def\Ptest{ }\def\Pcnt{ }\def\Pstr{171--196}%
\def\Dtest{ }\def\Dstr{1991}%
\def\Qtest{ }\def\Qstr{access via "jaffe sections commutative"}%
\def\Htest{ }\def\Hstr{4}%
\Refformat\egroup%

\bgroup\Resetstrings%
\def\Loccittest{}\def\Abbtest{ }\def\Capssmallcapstest{}\def\Edabbtest{
}\def\Edcapsmallcapstest{}\def\Underlinetest{ }%
\def\NoArev{1}\def\NoErev{0}\def\Acnt{1}\def\Ecnt{0}\def\acnt{0}\def\ecnt{0}%
\def\Ftest{ }\def\Fstr{9}%
\def\Atest{ }\def\Astr{Mitchell\Revcomma B\Initper }%
\def\Ttest{ }\def\Tstr{Theory of Categories}%
\def\Itest{ }\def\Istr{Academic Press}%
\def\Ctest{ }\def\Cstr{New York}%
\def\Dtest{ }\def\Dstr{1965}%
\def\Qtest{ }\def\Qstr{access via "mitchell categories"}%
\Refformat\egroup%

\bgroup\Resetstrings%
\def\Loccittest{}\def\Abbtest{ }\def\Capssmallcapstest{}\def\Edabbtest{
}\def\Edcapsmallcapstest{}\def\Underlinetest{ }%
\def\NoArev{1}\def\NoErev{0}\def\Acnt{1}\def\Ecnt{0}\def\acnt{0}\def\ecnt{0}%
\def\Ftest{ }\def\Fstr{10}%
\def\Atest{ }\def\Astr{Raynaud\Revcomma M\Initper }%
\def\Ttest{ }\def\Tstr{{\rm Groupes algebriques unipotents.  Extensions entre
groupes unipotents et groupes de type multiplicatif, \EXPOSE\ XVII in {\it\SGA}
(SGA 3)}}%
\def\Stest{ }\def\Sstr{Lecture Notes in Mathematics}%
\def\Itest{ }\def\Istr{Spring\-er-Ver\-lag}%
\def\Ctest{ }\def\Cstr{New York}%
\def\Vtest{ }\def\Vstr{152}%
\def\Dtest{ }\def\Dstr{1970}%
\def\Ptest{ }\def\Pcnt{ }\def\Pstr{532--631}%
\def\Qtest{ }\def\Qstr{access via "raynaud unipotents multiplicatif" (was
raynaud unipotents)}%
\Refformat\egroup%

\bgroup\Resetstrings%
\def\Loccittest{}\def\Abbtest{ }\def\Capssmallcapstest{}\def\Edabbtest{
}\def\Edcapsmallcapstest{}\def\Underlinetest{ }%
\def\NoArev{1}\def\NoErev{0}\def\Acnt{1}\def\Ecnt{0}\def\acnt{0}\def\ecnt{0}%
\def\Ftest{ }\def\Fstr{11}%
\def\Atest{ }\def\Astr{Rosenlicht\Revcomma M\Initper }%
\def\Ttest{ }\def\Tstr{Some basic theorems on algebraic groups}%
\def\Jtest{ }\def\Jstr{Amer. J. Math.}%
\def\Vtest{ }\def\Vstr{78}%
\def\Dtest{ }\def\Dstr{1956}%
\def\Ptest{ }\def\Pcnt{ }\def\Pstr{401--443}%
\def\Qtest{ }\def\Qstr{access via "rosenlicht basic"}%
\Refformat\egroup%

\bgroup\Resetstrings%
\def\Loccittest{}\def\Abbtest{ }\def\Capssmallcapstest{}\def\Edabbtest{
}\def\Edcapsmallcapstest{}\def\Underlinetest{ }%
\def\NoArev{1}\def\NoErev{0}\def\Acnt{1}\def\Ecnt{0}\def\acnt{0}\def\ecnt{0}%
\def\Ftest{ }\def\Fstr{12}%
\def\Atest{ }\def\Astr{Serre\Revcomma J-P\Initper }%
\def\Ttest{ }\def\Tstr{Algebraic Groups and Class Fields}%
\def\Itest{ }\def\Istr{Spring\-er-Ver\-lag}%
\def\Ctest{ }\def\Cstr{New York}%
\def\Dtest{ }\def\Dstr{1988}%
\def\Qtest{ }\def\Qstr{access via "serre algebraic groups and class fields"}%
\def\Xtest{ }\def\Xstr{descent: p. 102-104}%
\Refformat\egroup%

\bgroup\Resetstrings%
\def\Loccittest{}\def\Abbtest{ }\def\Capssmallcapstest{}\def\Edabbtest{
}\def\Edcapsmallcapstest{}\def\Underlinetest{ }%
\def\NoArev{1}\def\NoErev{0}\def\Acnt{1}\def\Ecnt{0}\def\acnt{0}\def\ecnt{0}%
\def\Ftest{ }\def\Fstr{13}%
\def\Atest{ }\def\Astr{Waterhouse\Revcomma W\Initper \Initgap C\Initper }%
\def\Ttest{ }\def\Tstr{Introduction to Affine Group Schemes}%
\def\Itest{ }\def\Istr{Spring\-er-Ver\-lag}%
\def\Ctest{ }\def\Cstr{New York}%
\def\Dtest{ }\def\Dstr{1979}%
\def\Qtest{ }\def\Qstr{access via "waterhouse introduction"}%
\Refformat\egroup%

\end{document}